
\def\LaTeX{%
  \let\Begin\begin
  \let\End\end
  \let\salta\relax
  \let\finqui\relax
  \let\futuro\relax}

\def\UK{\def\our{our}\let\sz s}
\def\USA{\def\our{or}\let\sz z}

\UK 



\LaTeX

\USA


\salta

\documentclass[twoside,12pt]{article}
\setlength{\textheight}{24cm}
\setlength{\textwidth}{16cm}
\setlength{\oddsidemargin}{2mm}
\setlength{\evensidemargin}{2mm}
\setlength{\topmargin}{-15mm}
\parskip2mm


\usepackage[usenames,dvipsnames]{color}
\usepackage{amsmath}
\usepackage{amsthm}
\usepackage{amssymb,bbm}
\usepackage[mathcal]{euscript}

\usepackage{cite}
\usepackage{hyperref}
\usepackage[shortlabels]{enumitem}

\usepackage[ulem=normalem,draft]{changes}
%
%

%
 
\definecolor{ciclamino}{rgb}{0.5,0,0.5}
\definecolor{blu}{rgb}{0,0,0.7}
\definecolor{rosso}{rgb}{0.85,0,0}
\def\jurgen #1{{\color{green}#1}}
\def\juerg #1{{\color{black}#1}}
\def\an #1{{\color{blue}#1}}
\def\last #1{{\color{rosso}#1}}
\def\pier #1{{\color{red}#1}}

\def\jurgen #1{#1}
\def\juerg #1{#1}
\def\pier #1{#1}
\def\an #1{#1}
\def\last #1{#1}




\bibliographystyle{plain}


%

\finqui

\def\Beq{\Begin{equation}}
\def\Eeq{\End{equation}}

\def\Bthm{\Begin{theorem}}
\def\Ethm{\End{theorem}}

\def\Brem{\Begin{remark}\rm}
\def\Erem{\End{remark}}

\def\Bdim{\Begin{proof}}
\def\Edim{\End{proof}}
\def\Bcenter{\Begin{center}}
\def\Ecenter{\End{center}}
\let\non\nonumber




\def\step #1 \par{\medskip\noindent{\bf #1.}\quad}
\def\jstep #1: \par {\vspace{2mm}\noindent\underline{\sc #1 :}\par\nobreak\vspace{1mm}\noindent}

\def\aand{\quad\hbox{and}\quad}
\def\Lip{Lip\-schitz}
\def\Holder{H\"older}
\def\Frechet{Fr\'echet}
\def\Poinc{Poincar\'e\pier{{}--Wirtinger{}}}
\def\Lady{Lady\v zhenskaya}
\def\lhs{left-hand side}
\def\rhs{right-hand side}




\def\multibold #1{\def\arg{#1}%
  \ifx\arg\pto \let\next\relax
  \else
  \def\next{\expandafter
    \def\csname #1#1\endcsname{{\boldsymbol #1}}%
    \multibold}%
  \fi \next}

\def\pto{.}

\def\multical #1{\def\arg{#1}%
  \ifx\arg\pto \let\next\relax
  \else
  \def\next{\expandafter
    \def\csname cal#1\endcsname{{\cal #1}}%
    \multical}%
  \fi \next}

\def\multigrass #1{\def\arg{#1}%
  \ifx\arg\pto \let\next\relax
  \else
  \def\next{\expandafter
    \def\csname grass#1\endcsname{{\mathbb #1}}%
    \multigrass}%
  \fi \next}


\def\multimathop #1 {\def\arg{#1}%
  \ifx\arg\pto \let\next\relax
  \else
  \def\next{\expandafter
    \def\csname #1\endcsname{\mathop{\rm #1}\nolimits}%
    \multimathop}%
  \fi \next}

\multibold
qweryuiopasdfghjklzxcvbnmQWERTYUIOPASDFGHJKLZXCVBNM.  

\multical
QWERTYUIOPASDFGHJKLZXCVBNM.

\multigrass
QWERTYUIOPASDFGHJKLZXCVBNM.

\multimathop
diag dist div dom mean meas sign supp .


\def\accorpa #1#2{\eqref{#1}--\eqref{#2}}
\def\Accorpa #1#2 #3 {\gdef #1{\eqref{#2}--\eqref{#3}}%
  \wlog{}\wlog{\string #1 -> #2 - #3}\wlog{}}


\def\separa{\noalign{\allowbreak}}

\def\somma #1#2#3{\sum_{#1=#2}^{#3}}

\def\graffe #1{\mathopen\{#1\mathclose\}}
\def\<#1>{\mathopen\langle #1\mathclose\rangle}
\def\norma #1{\mathopen \| #1\mathclose \|}

\def\checkmmode #1{\relax\ifmmode\hbox{#1}\else{#1}\fi}
\def\aeO{\checkmmode{a.e.\ in~$\Omega$}}
\def\aeG{\checkmmode{a.e.\ on~$\Gamma$}}
\def\aeQ{\checkmmode{a.e.\ in~$Q$}}
\def\aet{\checkmmode{a.e.\ in~$(0,T)$}}
\def\aat{\checkmmode{for a.a.\ $t\in(0,T)$}}

\let\hat\widehat
\def\cpto{\,\cdot\,}

\def\iot {\int_0^t}
\def\ioT {\int_0^T}
\def\intQt{\int_{Q_t}}
\def\intQ{\int_Q}
\def\iO{\int_\Omega}

\def\dt{\partial_t}
\def\dn{\partial_{\nn}}
\def\ddt{\frac d{dt}}

\def\0{{\boldsymbol {0} }}

\let\emb\hookrightarrow


\let\erre\grassR




\def\genspazio #1#2#3#4#5{#1^{#2}(#5,#4;#3)}
\def\spazio #1#2#3{\genspazio {#1}{#2}{#3}T0}

\def\L {\spazio L}
\def\H {\spazio H}
\def\W {\spazio W}

\def\C #1#2{C^{#1}([0,T];#2)}


\def\Lx #1{L^{#1}(\Omega)}
\def\Hx #1{H^{#1}(\Omega)}
\def\Wx #1{W^{#1}(\Omega)}

\def\Cx #1{C^{#1}(\overline\Omega)}

\def\LQ #1{L^{#1}(Q)}

\def\CQ #1{C^{#1}(\overline Q)}

\def\Luno{\Lx 1}
\def\Ldue{\Lx 2}

\def\Huno{\Hx 1}
\def\Hdue{\Hx 2}



\let\theta\vartheta

\let\ph\phi
\let\phi\varphi

\let\al\alpha
\let\sig\sigma
\let\lam\lambda

\let\TeXchi\chi                         
\newbox\chibox
\setbox0 \hbox{\mathsurround0pt $\TeXchi$}
\setbox\chibox \hbox{\raise\dp0 \box 0 }
\def\chi{\copy\chibox}



\def\normaV #1{\norma{#1}_V}
\def\normaW #1{\norma{#1}_W}

\def\phiz{\phi_0}
\def\az{a_0}
\def\nz{n_0}
\def\sigz{\sig_0}
\def\rz{r_0}
\def\deltaz{\delta_0}

\def\Vp{{V^*}}

\def\soluz{(\phi,\mu,a,n,\sig)}
\def\soluzstar{(\phi^*,\mu^*,a^*,n^*,\sig^*)}
\def\soluzk{(\phik,\muk,\ak,\nk,\sigk)}
\def\soluzl{(\psi,\eta,\al,\nu,\omega)}
\def\soluzh{(\phih,\muh,\ah,\nh,\sigh)}
\def\soluzF{(\ph,\rho,\gamma,\lam,\an{\xi})}
\def\soluzad{(\p1,\p2,\p3,\p4,\p5)}

\def\vbar{\overline v}
\def\phibar{\overline\phi}
\def\mubar{\overline\mu}

\def\phizbar{\overline\phiz}
\def\gbar{\overline g}

\def\rbar{r'}
\def\rbarmin{\rz - (\phizbar-\rz)^-\!-R}
\def\rbarmax{\rz + (\phizbar-\rz)^+\!+R}

\def\h{{\mathbbm h}}
\def\chiphi{\chi_\phi}
\def\chia{\chi_a}
\def\cphi{c_\phi}
\def\cn{c_n}
\def\csig{c_\sig}
\def\Sn{S_n}
\def\Ssig{S_\sig}
\def\Funo{F_1}
\def\Fdue{F_2}

\def\ustar{u^*}
\def\phistar{\phi^*}
\def\mustar{\mu^*}
\def\astar{a^*}
\def\nstar{n^*}
\def\sigstar{\sig^*}

\def\uk{u_k}
\def\phik{\phi_k}
\def\muk{\mu_k}
\def\ak{a_k}
\def\nk{n_k}
\def\sigk{\sig_k}

\def\phih{\phi^h}
\def\muh{\mu^h}
\def\ah{a^h}
\def\nh{n^h}
\def\sigh{\sig^h}

\def\umax{u_{\rm max}}

\def\J{{\cal J}}

\def\Uad{{\cal U}_{\rm ad}}
\def\phiQ{\phi_Q}
\def\phiO{\phi_\Omega}

\def\b #1{b_{#1}}
\def\p #1{p_{#1}}
\def\f #1#2{f_{#1,#2}}

\def\CO{C_\Omega}


\usepackage{amsmath}
\DeclareFontFamily{U}{mathc}{}
\DeclareFontShape{U}{mathc}{m}{it}%
{<->s*[1.03] mathc10}{}

\DeclareMathAlphabet{\mathscr}{U}{mathc}{m}{it}

\newcommand{\ov}[1]{\overline{#1}}

\Begin{document}


%
\title{\pier{Solvability and} optimal control of a multi-species Cahn--Hilliard--Keller--Segel tumor growth model}
\author{}
\date{}
\maketitle
\Bcenter
\vskip-1.5cm
{\large\sc Pierluigi Colli$^{(1)}$}\\
{\normalsize e-mail: {\tt pierluigi.colli@unipv.it}}\\[0.25cm]
{\large\sc Gianni Gilardi $^{(1)}$}\\
{\normalsize e-mail: {\tt gianni.gilardi@unipv.it}}\\[0.25cm]
{\large\sc Andrea Signori$^{(2)}$}\\
{\normalsize e-mail: {\tt andrea.signori@polimi.it}}\\[0.25cm]
{\large\sc J\"urgen Sprekels$^{(3)}$}\\
{\normalsize e-mail: {\tt juergen.sprekels@wias-berlin.de}}\\[.5cm]
$^{(1)}$
{\small Dipartimento di Matematica ``F. Casorati'', Universit\`a di Pavia}\\
{\small and Research Associate at the IMATI -- C.N.R. Pavia}\\
{\small via Ferrata 5, I-27100 Pavia, Italy}\\[.3cm] 
$^{(2)}$
{\small Dipartimento di Matematica, Politecnico di Milano}\\
{\small via E. Bonardi 9, I-20133 Milano, Italy}\\
{\small \pier{and} \an{Alexander von Humboldt Research Fellow}}
\\[.3cm] 
$^{(3)}$
{\small Department of Mathematics}\\
{\small Humboldt-Universit\"at zu Berlin}\\
{\small Unter den Linden 6, D-10099 Berlin, Germany}\\
{\small and}\\
{\small Weierstrass Institute for Applied Analysis and Stochastics}\\
{\small Mohrenstrasse 39, D-10117 Berlin, Germany}\\[10mm]

\Ecenter
\Begin{abstract}
\noindent 
This paper investigates an optimal control problem associated with a two-dimen\-sional 
multi-species Cahn\last{--}Hilliard\last{--}Keller\last{--}Segel tumor growth model, which incorporates complex biological processes such as species diffusion, chemotaxis, angiogenesis, and nutrient consumption, resulting in a highly 
nonlinear system of nonlinear partial differential equations.
The modeling derivation and corresponding analysis have been addressed in a previous contribution. Building on this foundation, the scope of this study involves investigating a distributed control problem with the goal of optimizing a tracking-type cost 
functional. This latter  aims to minimize the deviation of tumor cell location from desired target 
configurations while penalizing the costs associated with implementing control measures, akin to introducing a suitable medication.
Under appropriate mathematical assumptions, we demonstrate that sufficiently regular solutions exhibit continuous dependence on the control variable. Furthermore, we establish the existence of optimal controls and characterize the first-order necessary optimality conditions through a suitable variational inequality.


\vskip3mm
\noindent 
{\bf AMS (MOS) Subject Classification:} {
		35K55, 
        \pier{35K52}, 
		49J20, 
		49K20, 
		49K40. 
		}
\End{abstract}

\pagestyle{myheadings}
\newcommand\testopari{\sc Colli -- Gilardi -- Signori -- Sprekels}
\newcommand\testodispari{\sc Optimal control of a Cahn--Hilliard--Keller--Segel model}
\markboth{\testopari}{\testodispari}
%

\section{Introduction}
\label{INTRO}
\setcounter{equation}{0}

This paper investigates an optimal control problem associated with a multi-species Cahn--Hilliard--Keller--Segel tumor growth model in a two-dimensional spatial domain $\Omega \subset \mathbb{R}^2$ over a given final time $T > 0$. 
The problem we aim to analyze consists of a distributed optimal control problem associated with an initial-boundary value problem:
\begin{alignat}{2}
  & \dt\phi - \Delta\mu + \chiphi \Delta n
  = - m\phi + \h(\phi)
  \quad && \text{in $Q:=\Omega \times (0,T)$},
  \label{Iprima}
  \\
  & {- \Delta\phi} + F'(\phi) 
  = \mu
  \quad && \text{in $Q$},
  \label{Iseconda}
  \\
  & \dt a - \Delta a + \chia \div (a \nabla\sig)
  = a - a^2 + u
  \quad && \text{in $Q$},
  \label{Iterza}
  \\
  & \dt n - \Delta{n} - \pier{\chiphi n}
  = \pier{\Sn}
  \quad && \text{in $Q$},
  \label{Iquarta}
  \\
  & \dt\sig - \Delta{\sig} - \pier{\chia a}
  = \pier{\Ssig}
  \quad && \text{in $Q$},
  \label{Iquinta}
  \\
  & \dn\phi
  = \dn\mu
  = \dn a
  = \dn n
  = \dn\sig
  = 0 
  \quad && \text{on $\Sigma:=\partial\Omega \times (0,T)$},
  \label{Ibc}
  \\
  & \phi(0) = \phiz, 
  \quad 
  a(0) = \az,
  \quad 
  n(0) = \nz,
  \quad 
  \sig(0) = \sigz
  \quad && \text{in $\Omega.$}
  \label{Icauchy}
\end{alignat}
\Accorpa\Ipbl Iprima Icauchy
The primary variables in the system are \( \phi, \mu, a, n, \) and \( \sigma \). These represent the density of tumor cells \( \phi \), the chemical potential \( \mu \), an angiogenetic phase composed of tumor-induced new vasculature with volume fraction \( a \), a nutrient or signaling molecule \( n \), and a concentration \( \sigma \) affecting tumor growth dynamics. The positive constants \( \chi_\phi \) and \( \chi_a \) represent chemotaxis parameters quantifying the sensitivity of biological entities to chemical gradients. In the second equation \eqref{Iseconda}, \( F' \) denotes the derivative of a configuration potential \( F \) characterized by a double-well shape. Prototypical choices for this latter include the so-called \textit{classical regular potential} and the \textit{logarithmic potential} defined as follows:
\begin{align}
    F_{reg}(r) & := \pier{\frac{c_1}4 r^2 (r - 1)^2}, \quad r \in \mathbb{R},
    \label{regpot} \\
    F_{log}(r) & := r \ln r + (1 - r) \ln (1 - r) + \pier{c_2} \, r (1 - r), \quad r \in (0,1),
    \label{logpot}
\end{align}
where \pier{\( c_1, \, c_2  \) are two positive real coefficients}. The mass of the tumor, represented by $\phi$, is not conserved, as indicated by the presence of a source term \( -m\phi + \h(\phi) \) on the right-hand side of the first equation, where \( \h \) represents a smooth real function and \( m \) is a positive constant. Chemotaxis is modeled via a Keller--Segel type (cf.,~e.g.,\pier{\cite{KS}}) coupling, specifically through the nonlinear term \( \chia \div(a \nabla \sigma) \) in the third equation. The logistic source term \( a - a^2 \) in \eqref{Iterza} for the nutrient variable \( a \) is a common choice in Keller--Segel models to prevent solution blow-up in finite times, see, e.g., \cite{DP, W_blow, HH, RSchS, GSigSpr} and the references therein. Finally, \( \phiz, \az, \nz, \) and \( \sigz \) denote prescribed initial data for these variables, whereas \pier{$S_n$ and $ S_\sig$} stand for suitable source terms \last{depending on the solution variables}, details of which will be provided later on.

The model \Ipbl\ originates from variational principles and was introduced in \cite{Agosti}\pier{.
The} postulated free energy of the system, which is defined as the internal energy minus the entropy, is given by
\begin{align}
  \calE(\phi,a,n,\sig)
  & = \iO a (\ln a - 1)
  - \chiphi \iO n \phi
  - \chia \iO a \sig
  \non
  \\
  & \quad + \frac 12 \iO |\nabla n|^2
  + \frac 12 \iO |\nabla\sig|^2
  + \frac 12 \iO |\nabla\phi|^2
  + \iO F(\phi).
  \label{Ienergy}
\end{align}
In \cite{Agosti}, the modeling derivation and numerical simulations (cf. \an{\cite{AZAMM}}) aim to optimize model parameters \pier{and support clinical decision-making,} whereas the corresponding mathematical analysis is addressed in \cite{AS}, where the existence of weak solutions was shown in two and three dimensions, and regularity results and uniqueness of regular enough solutions were proved in the two-dimensional setting. Here, we aim at considering an optimal control problem, where the distributed control $u$ enters in the form of a source term in \eqref{Iterza}.
The minimization problem we want to study consists in minimizing a suitable  {\it cost functional} that we postulate to be of tracking type form\pier{, and expressed by}
\Beq
	\J (\phi, u)
	:= \frac {\b1}2 \intQ |\phi - \phiQ|^2
	+ \frac {\b2}2  \iO |\phi(T) - \phiO|^2
	+ \frac {\b3}2 \intQ |u|^2 \,,
	\label{Icost}
\Eeq
where the coefficients $\b i$ are given nonnegative numbers\pier{, with $b_3>0$,} and $\phiQ$ and $\phiO$ are given functions on $Q$ and~$\Omega$, respectively, representing clinical targets.
Besides, we \jurgen{constrain the} control variables to belong to the set of {\it admissible controls} defined by
\begin{align}
	\Uad : = \big\{ u \in \calU := \LQ\infty : 0 \leq u \leq u_{\rm max} \ \aeQ \big\} \,,
	\label{IUad}
\end{align}
where $\umax\in L^\infty(Q)$ is a prescribed {nonnegative function}.
Then the control problem {under investigation} can be formulated as {follows}:
\begin{align}
  & \hbox{Minimize} \ \ \J (\phi, u) \ \ \text{subject to {$u\in\Uad$ and to the constraint that}}
  \non
  \\
  & \quad \hbox{{$\soluz$ is the solution to the system \Ipbl}}.
  \label{Icontrol}
\end{align}   
Tumor growth models based on the phase field approach have gained significant popularity. While not exhaustive, we refer interested readers to \cite{CGH, CL, GLS, GLSS} and the references therein. Several studies within this framework consider the influence of velocity effects on the mixture dynamics, utilizing Darcy's law and the Brinkman equation. For detailed discussions on these topics, see \cite{CGSS1, EG, KS2, GLSS}.
The incorporation of chemotaxis, particularly through the Keller--Segel coupling, represents a relatively recent advancement in phase field models. This coupling has been explored in studies such as \cite{RSchS, Agosti, GSigSpr}.
Finally, regarding the optimal control problem, we refer readers to \cite{S, CGRS-oc, GLS_OPT, CSS, CSS2, KE, S_time} for comprehensive discussions and analyses.

\medskip
The plan of the paper is as follows.
In the next section, we state the problem in a precise form
and present our results.
The existence and the uniqueness of the solution to the state system,
as well as proper stability and continuous dependence estimates,
are proved in Sections \ref{EXISTENCE} and~\ref{UNIQUENESS}.
A~technical result is given in Section~\ref{AUXILIARY}\last{, preparing} the study of the control problem made in the last two sections,
where we prove the existence of an optimal control and
we establish first\last{-}order necessary conditions for optimality
in terms of the solution to the adjoint problem.


\section{Statement of the problem and results}
\label{STATEMENT}
\setcounter{equation}{0}

Throughout the paper, $\Omega$ is a \juerg{bounded open subset of $\erre^2$ having a smooth boundary $\Gamma:=\partial\Omega$}.
The symbols $|\Omega|$ and $\dn$ \last{denote} the measure of $\Omega$
and the derivative in the direction of the outward unit normal vector $\nn$ 
on $\Gamma$\last{, respectively}. \juerg{With a prescribed final time $T>0$, we set} 
\Beq
  Q := \Omega\times(0,T)\juerg{, \quad \Sigma:=\Gamma\times(0,T).} 
  \label{defQ}
\Eeq
Given a Banach space~$X$, we denote by $\norma\cpto_X$ 
both its norm and the norm in any power of~$X$,
with the exceptions of the space $H$ introduced below and of the Lebesgue spaces $\Lx p$ ($1\leq p\leq+\infty$),
for which we use the symbol $\norma\cpto_p$. Sometimes, this symbol also denotes the norm in~$\LQ p$.
Moreover, \juerg{in order to simplify the notation}, we still write~$X$ (i.e.,~we avoid the exponent) when dealing with some power of~$X$.
Then, we introduce the shorthands
\Beq
  H := \Ldue, \quad
  V := \Huno
  \aand
  W := \graffe{ v\in\Hdue: \dn v=0 \ \aeG},
  \label{defspazi}
\Eeq
and endow these spaces with their natural norms.
For simplicity, we write $\norma\cpto$ instead of~$\norma\cpto_H$.
Moreover, we denote by $\Vp$ and $\<\cpto,\cpto>$ the dual space of $V$ and the duality pairing between $\Vp$ and~$V$, respectively,
and we identify $H$ with a subspace of $\Vp$ in the usual way, i.e., \juerg{such} that
$\,\<w,v>=\iO wv\,$ for every $w\in H$ and $v\in V$.
Hence, we have the continuous, dense, and compact embeddings
\Beq
  V \emb H \emb V^*\,, 
  \non
\Eeq
yielding that $(V,H,\Vp)$ is a Hilbert triple.

\medskip

At this point, we are ready to introduce our assumptions on the structure of the state system \an{which} \pier{involve}, in particular,  \pier{a specific} choice of the functions $\pier{\Sn}$ and $\pier{S_\sigma}$ 
that appear in equations \eqref{Iquarta} and~\eqref{Iquinta}.
We assume~that
\begin{align}
  & m \in (0,+\infty) , \quad \pier{\chiphi, \chia \in (0,1),} \quad
  \cphi , \cn , \csig , c_0 \in \erre,
  \label{HPconst}
  \\
  & \hbox{$\h:\erre\to\erre$  is \an{such that $\h \in W^{2,\infty}(\erre)$},}
  \label{HPh}
  \\
  &  \pier{{}\Sn : = \cphi \phi + \cn n + \csig \sig + c_0, \quad  
  \Ssig :=  1-\sig - a \,\sig  }.
  \label{HPS}
\end{align}
\Accorpa\HPstructure HPconst HPS
As for the potential \pier{$F$}, we confine ourselves to \pier{some conditions}
that generalize the cases 
of the classical and logarithmic potentials \last{\eqref{regpot} and \eqref{logpot}}.
In particular, we \pier{ignore the possibility of extending the latter to $[0,1]$ by continuity
and prescribe suitable regularities on $F$ in the open interval that is taken as the domain $D(F)$ 
in any case and is actually the effective domain of the derivative $F'$.}
We~assume: 
\begin{align}
  & \hbox{Either \ $D(F)=\erre$ \ or \ $D(F)=(0,1)$},
  \,\juerg{\mbox{and it holds that}}
  \non
  \\
  &\quad{} F = \Funo + \Fdue \,\,\hbox{ with \juerg{functions} } \,\,
  \Funo\,,\Fdue : D(F) \to \erre \,\, \hbox{ of class $C^4$,}
  \non
  \\
  & \quad{}\hbox{\juerg{where} $\,\Funo\,$ is convex and $\,\Fdue'\,$ is \Lip\ continuous}.
  \label{HPsplit}
  \\[1mm]
  & \hbox{\pier{If \ $D(F)=\erre$, \ then $F$ satisfies}} \quad
  \lim_{|r|\to+\infty} r^{-2} F(r) = +\infty\,;
  \label{HPregpot} 
  \\
  &\pier{\hbox{if \ $D(F)=(0,1)$, \ then}
  \quad
  \lim_{r\searrow  0} F'(r) = -\infty ,
  \quad
  \lim_{r \nearrow 1} F'(r) = +\infty , \quad \hbox{and there is}}
  \non
  \\
  &\quad{}\pier{\hbox{a constant $C_F$ such that } \, |F_1'' (r)|\leq e^{C_F (|F_1'(r)| +1)}\, 
  \hbox{ for all } r\in (0,1).}
  \label{HPlogpot}
 \end{align}
\Accorpa\HPF HPsplit HPlogpot
\Accorpa\HPall HPconst HPlogpot
We notice that these assumptions imply that $F$ is bounded from below \pier{and also ensure the existence of \juerg{some} $\rz \in\erre$ satisfying}
\Beq
  \pier{\rz \in D(F)
  \aand
  \Funo'(\rz) = 0 \,.}
  \label{defrz}
\Eeq
\pier{As the reader can directly check\last{,} it turns out that the growth condition in \eqref{HPlogpot} is satisfied by the convex part of the logarithmic potential in \eqref{logpot}.} 

\juerg{For the control variable $u$, we assume that
\Beq
  u \in \LQ\infty
  \quad \hbox{satisfies} \quad
  0 \leq u \leq\umax \quad \aeQ \,,
  \label{HPu}
\Eeq
\Accorpa\HPu HPu HPumax
where}
\Beq
 \juerg{ \umax \in \LQ\infty
  \quad \hbox{is nonnegative.} }
  \label{HPumax}
\Eeq

To introduce our assumptions on the initial data,
we use the following general notation for the mean value: we~set
\Beq
  \vbar := \frac 1{|\Omega|} \, \iO v
  \quad \hbox{for $v\in\Luno$} \,.
  \label{mean}
\Eeq
The same symbol will be used in the sequel even for time-dependent functions.
Then, denoting by  $(\cpto)^\pm$ the positive and negative parts, we assume~that
\begin{align}
  & \pier{\phiz \in W  \, \hbox{ with range in $D(F)$; \ moreover, \,  
   $\phiz \in \Hx3$ \, if \, $D(F) = \erre $,}}
  \non
  \\
  & \quad{} \pier{\phiz \in \Hx4 \, \hbox{ and }\, \mu_0:=\last{-\Delta \phi_0} + F'(\phi_0) \in W \hbox{ \, if \, $D(F) = (0,1) $.}}
  \label{HPphiz}
  \\[1mm]
  & \pier{\hbox{$\rbarmin$ \ and \ $\rbarmax$ \ belong to $D(F)$},}
  \non \\
  &\quad \pier{\hbox{where} \quad
  R := \frac 1 m  \sup_{r\in \erre} |\h (r) - m\rz | \,.}
  \label{HPinterior}
  \\[1mm]
  & \az \in \an{V} 
  \aand
  \az > 0 \quad \aeO \,.
  \label{HPaz}
  \\
  & \nz , \sigz \in W 
  \aand
  0 \leq \sigz \leq 1 \quad \hbox{in $\Omega$} \,.
  \label{HPsigz}
\end{align}
\Accorpa\HPz HPphiz HPsigz
Of course, the \pier{assumption in \eqref{HPinterior} yields a restriction only in the case when $D(F)$ is bounded.}

Finally, we are in a position to introduce our formulation of the state system.
Even though some of the equations could be written in the strong form used in the Introduction,
we prefer to present the whole problem in a variational \an{framework}.
We look for a quintuple $\soluz$ with the properties
\pier{\begin{align}
  & \phi \in \calY_1 := \H1 {\pier{V}} \cap \L\infty {W\cap\Hx3} \cap \CQ0
  \non
  \\
  &\quad {}\hbox{and }\ \phi \in D(F) \ \ \aeQ \,,
  \label{regphi}
  \\[1mm]
  & \mu \in \calY_2 := \L\infty {\pier{V}}\cap \L2 {W\cap\Hx3} \,,
  \label{regmu}
  \\[1mm]
  & a \in \calY_3 := 
   \H1H \cap \L\infty V \cap \L2W\,
  \non
  \\
  & \aand
  a \pier{{}> 0{}} \quad \aeQ \,,
  \label{rega}
  \\[1mm]
  & n \in \calY_4 := \H1{V}\cap\L\infty W\cap\L2{\Hx{3}} \cap \last{\calZ}
 \non
  \\
  &\quad\hbox{where }\
  \last{\calZ} :=\W{1,4}{\Lx4}\cap\L4{\Wx{2,4}}\,,
  \label{regn}
  \\[1mm]
  & \sig \in \calY_4 \aand
  0 \leq \sig \leq 1 \quad \aeQ \,,
  \label{regsig}
\end{align}}%
\Accorpa\Regsoluz regphi regsig
that solves the variational equations
\begin{align}
&\iO \dt\phi \,v +\iO \nabla\mu\cdot\nabla v
  - \chiphi \iO \nabla n \cdot \nabla v
  = - m \iO \phi v
  + \iO \h(\phi) \, v \,,
  \label{prima}
  \\
  & \iO \nabla\phi \cdot \nabla v
  + \iO F'(\phi) \, v
  = \iO \mu v \,,
  \label{seconda}
  \\
  & \pier{\iO \dt a \, v}
  + \iO \nabla a \cdot \nabla v
  - \chia \iO a \nabla\sig \cdot \nabla v
  = \iO \bigl( a - a^2 + u \bigr) v \,,
  \label{terza}
  \\
  & \iO \dt n \, v
  + \iO \nabla n \cdot \nabla v
  - \chiphi \iO \phi v
  = \iO \Sn \, v 
  \non
  \\
  & \quad \hbox{where} \quad 
  \Sn = \cphi \phi + \cn n + \csig \sig + c_0 \,,
  \label{quarta}
  \\
  & \iO \dt\sig \, v
  + \iO \nabla\sig \cdot \nabla v
  = \iO \bigl( (1-\sig) + a (\chia-\sig) \bigr) v \,,
  \label{quinta}
\end{align}
\juerg{for} every $v\in V$ and \aet,
and satisfies the initial condition
\Beq
  (\phi,a,n,\sig)(0) = (\phiz,\az,\nz,\sigz) 
  \quad \aeO \,.
  \label{cauchy}
\Eeq
\Accorpa\Pbl prima cauchy

\Brem
\label{Remregularity}
We notice that the regularity properties \eqref{rega} and \eqref{regsig} imply
that both $a$ and $\nabla\sig$ are $L^4$ functions (see the forthcoming \eqref{embeddings}),
so that all \juerg{of the terms occurring} in \eqref{terza} are meaningful.
\juerg{We also point out that by virtue of \Regsoluz\ all of} the above equations may be written in their strong form.
\pier{From \eqref{regphi}} we \juerg{also have} that
\begin{align}
  &\phi \in \L\infty{\Hx3}\,,
  \quad \hbox{whence} \quad
  \nabla\phi \in \LQ\infty\,,
  \label{nablaphibdd}
\end{align}
\pier{with} the corresponding norms \pier{that} are estimated by a constant like~$K_1$ \juerg{below}.
\Erem

\an{For convenience, we set the state space as} \pier{(cf.~\Regsoluz)}
\Beq 
\calY:= \calY_1\times\calY_2\times\calY_3\times\calY_4\times\pier{{}\calY_4}.
\label{regtot}
\Eeq
\juerg{Our first result} regards well-posedness and stability \last{of system \Ipbl}.

\Bthm
\label{Wellposedness}
Assume \pier{\HPall} on the structure, 
and \HPu\ and \HPz\ on the data.
Then there exists a unique quintuple $\soluz$ that satisfies \Regsoluz\ and solves problem \Pbl.
Moreover, the stability estimate and separation property
\begin{align}
  & \norma\soluz_{\an{\calY}}
  \leq K_1 \,,
  \label{stability}
  \\
  & r_- \leq \phi \leq r_+ \quad \aeQ \,,
  \label{separation}
\end{align}
hold true with constants $K_1>0$ and $r_\pm\in D(F)$
that depend only on $\Omega$, $T$, the structure of the system, the initial data, and~$\umax$.
In particular, \last{they are independent of $u$.}
\Ethm

Next, we have the following continuous dependence result.

\Bthm
\label{Contdep}
\juerg{Suppose} the assumptions of Theorem \ref{Wellposedness} regarding the structure and the initial data \juerg{are fulfilled, and}
let $u_i$, $i=1,2$, satisfy \eqref{HPu} and $(\phi_i,\mu_i,a_i,n_i,\sig_i)\an{\in {\calY}}$ be the corresponding solution.
Then the inequality
\begin{align}
  & \norma{\phi_1-\phi_2}_{\H1\Vp\cap\pier{{}\L\infty V\cap {}}\L2W}
  + \norma{\mu_1-\mu_2}_{\L2V}
  \non
  \\
  & {} \quad
  + \norma{a_1-a_2}_{\H1\Vp \an{\cap }\L\infty H\cap\L2V}
  + \norma{n_1-n_2}_{\H1H\pier{{}\cap\L\infty V{}}\cap\L2W}
  \non
  \\
  & {} \quad
  + \norma{\sig_1-\sig_2}_{\H1H\cap\L\infty V\cap\L2W}
  \leq K_2 \, \norma{u_1-u_2}_{\L2H}
  \label{contdep}
\end{align}
holds true with a constant $K_2>0$ 
that depends only on $\Omega$, $T$, the structure of the system, the initial data, and~$\umax$.
\Ethm

Once well-posedness is established, we can deal with the control problem presented in the Introduction
(see \accorpa{Icost}{Icontrol}).
\an{Let us} refer to the last two sections for \an{the} precise statements.
\an{Here, we  just mention} that we first prove the existence of an optimal control, i.e., of an element~$\ustar\in\Uad$
that satisfies
\Beq
  \calJ(\phistar,\ustar) \leq \calJ(\phi,u)
  \quad \hbox{for every $u\in\Uad$}
  \label{opt}
\Eeq
where $\phistar$ and $\phi$ are the first components of the solutions corresponding to $\ustar$ and~$u$, respectively.
Then, we derive a first-order necessary optimality condition for a given $\ustar\in\Uad$ to be an optimal control \an{in terms of a suitable variational inequality. Namely, $\ustar$ is an optimal control whether it fulfills}
\Beq
  \intQ (p_3 + \b3\ustar) (u-\ustar) \geq 0
  \quad \hbox{for every $u\in\Uad$}\,,
  \label{nc}
\Eeq
where $\p3$ is the third component of the solution to the adjoint problem introduced and discussed in Section~\ref{NECESSARY}.

\medskip

In performing our proofs, we often make use \pier{of \Holder's inequality, as well as of Young's inequality
\Beq
  yz \leq \frac \delta p |y|^p + \frac 1{p'} \delta^{- p'/p} \, |z|^{p'}
  \quad \hbox{for every $y,z\in\erre$ and $\delta>0$}\,,
  \label{young}
\Eeq
with $1< p,\,  p'<\infty$ conjugate exponents, i.e., $p+p' = p\, p'$.}  
Moreover, we recall the two-dimensional embeddings
\begin{align}
  & V \emb \Lx p \quad \hbox{for $p\in[1,+\infty)$}, \quad
  W \emb \Cx0, 
  \non
  \\
  & \mbox{and }\,  \L\infty H \cap \L2V \emb \LQ4 ,
  \label{embeddings}
\end{align}
and the corresponding inequalities
\begin{align}
  & \norma v_p \leq C_{\Omega,p} \, \normaV v
  \quad \hbox{for every $v\in V$}, \quad
  \norma v_\infty \leq \CO \, \normaW v
  \quad \hbox{for every $v\in W$}
  \non
  \\
  & \aand \norma v_{\LQ4} \leq C_{\Omega,T} \, \norma v_{\L\infty H\cap\L2V}
  \non
  \\
  & \quad \hbox{for every $v\in\L\infty H\cap\L2V$},
  \label{disug}
\end{align}
where $\CO$ depend only on $\Omega$ and $C_{\Omega,p}$ and $C_{\Omega,T}$ depend on $p$ and~$T$, in addition.
\pier{Furthermore, since the embeddings $V\subset H$ and $H\subset\Vp$ are compact,
we {obtain from Ehrling's lemma} the compactness inequality 
\Beq
  {\|v\|}
  \,\leq\, \delta\, \norma{\nabla v}
  + C_\delta \, \norma v_{\Vp}
  \quad \hbox{for every $v\in V$ and $\delta>0$},
  \label{compact}
\Eeq
with some $C_\delta>0 $ that depends only on~$\Omega$ and~$\delta$.}
We also account for the \Poinc\ inequality, \pier{an inequality from the elliptic regularity theory, and the two-di\-mensional \an{\Lady}\ inequality.} \an{Namely, we sometimes owe~to}
\begin{align}
  & \norma{v-\vbar} \leq \CO \, \norma{\nabla v}
  \quad \hbox{for every $v\in V$},
  \label{poincare}
  \\
  & \norma{\nabla v} \leq \CO \, (\norma v + \norma{\Delta v})
  \quad \hbox{for every $v\in W$},
  \label{elliptic}
  \\[1mm]
  & \norma v_4^2 \leq \CO \, \norma v \, \normaV v
  \quad \hbox{for every $v\in V$},
  \non
  \\
  &\quad \hbox{and } \ \norma{\nabla v}_4^2 \leq \CO \, \normaV v \, (\norma v + \norma{\Delta v})
  \quad \hbox{for every $v\in W$},
  \label{lady}
\end{align}
with the same constant $\CO$ as before, without loss of generality. \pier{We aim to point out that 
$v \mapsto \norma v + \norma{\Delta v}$ provides a norm in $W$ which is equivalent to the standard norm in $\Hx2$.}

We conclude this section by stating a convention that regards the constants appearing in the proofs of the \last{forthcoming} sections.
The small-case symbol $c$ denotes a generic constant
that depends only on the structure of the system, $\Omega$, $T$, the initial data,
and~$\umax$ (see~\eqref{HPumax}). 
In particular, the values of $c$ are independent of~$u$.
Notice that the meaning of $c$ may vary from line to line and even within the same line.
We use \an{capital letters} for precise constants we could refer~to.


\section{Existence and stability} 
\label{EXISTENCE}
\setcounter{equation}{0}

This section is devoted to the existence of a solution $\soluz$ to problem \Pbl\ 
that satisfies \pier{estimates~\eqref{stability} and~\eqref{separation}}.
\last{We mention that} \an{the well-posedness of} a \an{similar system} \pier{can be compared with the analysis developed} in \cite{AS}. 
\an{T}he system studied \an{there} has a \an{slightly} more general structure, 
but it does not contain the \an{control variable~$u$, acting as a} source term.
The authors \pier{of \cite{AS}} prove the existence of~a (unique regular) solution, satisfying 
the restrictions on the values of $a$ and $\sig$ given in our statement,
by means of a proper argument based on regularization, truncation and discretization
(see \cite[\pier{Theorems}~2.8--2.11]{AS}). \pier{On the other hand, here we construct 
our argumentation without giving the full detail of approximation
and just} \an{perform formal a priori estimates on \pier{the solution $\soluz$ to motivate 
the expected regularity and the stability estimates}. 
In particular, we \pier{point out the treatment of the new terms} involving the control $u$, which do not appear in the paper~\cite{AS}.}
\pier{In agreement with the specific form of the energy $\calE$ in~\eqref{Ienergy},} we assume \juerg{in the following} \an{the component $a$ to be positive.}

\step
\pier{Boundedness property}

We first prove that $\sig$ \juerg{attains its} values in~$[0,1]$.
To this end, we fix a monotone $C^1$ function $G:\erre\to\erre$ that grows linearly at infinity
and satisfies $G(r)<0$ for $r<0$, $G(r)=0$ for $r\in[0,1]$, and $G(r)>0$ for $r>1$,
and we test \eqref{quinta} by~$G(\sig)$.
If $\hat G$ is the \an{antiderivative} of $G$ that vanishes at zero, we obtain~that
\Beq
  \frac 12 \, \ddt \iO \hat G(\sig)
  + \iO G'(\sig) |\nabla\sig|^2
  = \iO \bigl( (1-\sig) + a(\chia-\sig) \bigr) G(\sig) 
  \quad \aet.
  \non
\Eeq
Since $a$ is nonnegative and $\chia\in(0,1)$, the \rhs\ is nonpositive.
\an{We then integrate over time and observe that}
our assumptions on~$\sigz$ (see~\eqref{HPsigz}) \juerg{imply} that $\hat G(\sigz)=0$.
Since both $\hat G$ and $G'$ are nonnegative, we conclude that $\hat G(\sig)$ vanishes identically\an{, entailing} that
\Beq
  0 \leq \sig \leq 1 \quad \aeQ .
  \label{stimasig}
\Eeq

\step
{\an{Control of the mean value}}

The next estimate regards the mean value of~$\phi$.
We test \eqref{prima} by the constant $1/|\Omega|$ and obtain~that
\Beq
  \ddt \, \phibar + m \, \phibar = \gbar\,,
  \quad \hbox{where} \quad
  g := \h(\phi) \,.
  \label{pier1}
\Eeq
\pier{In addition, recalling \eqref{defrz} we note that the constant function $v(t)=\rz$ fulfills the equation 
\Beq
  \ddt \, v + m\, v = m \,\rz \quad \hbox{in }\, (0,T)
  \label{pier2}
\Eeq
and the initial condition $v(0)=\rz$. Hence, taking the difference between 
\eqref{pier1} and \eqref{pier2}, and solving the resulting Cauchy problem for $\phibar - v$, we easily find that
\Beq
  \phibar(t) - \rz  = (\phizbar - \rz) \, e^{-mt} + \iot e^{-m(t-s)} \, (\gbar(s) - m \rz)  \, ds
  \quad \hbox{for every $t\in[0,T]$}.
  \non
\Eeq
Since $ \, |\gbar(s) - m \rz| \leq  \sup_{r\in \erre} |\h (r) - m\rz | = m R  $ \ for all $\, s\in [0,T] \, $
(see \eqref{HPinterior}), we easily conclude~that}
\Beq
  \pier{\rbarmin \leq \phibar(t) \leq \rbarmax
  \quad \hbox{for every $t\in[0,T]$}.}
  \label{meanok}
\Eeq

Finally, by recalling \eqref{HPinterior}, we claim that
there are positive constants $\deltaz$ and $C_0$ such~that
\begin{align}
  & \Funo'(r) (r-\rbar)
  \geq \deltaz \, |\Funo'(r)| - C_0
  \non
  \\[1mm]
  & \quad \hbox{for every $r\in \pier{D(F)}$ and $\rbar\in[\rbarmin,\rbarmax]$}.
  \label{trickMZ}
\end{align}
This is a generalization of the inequality proved in \cite[Appendix, Prop.~A.1]{MiZe}
in the case of a fixed $\rbar$. 
However, the proof also works in the present case with only minor changes 
since the values of $\rbar$ we are considering belong to a compact subset 
of the open interval~$D(F)$.

\medskip

At this point, we start performing the estimates in the direction of the expected regularity of the solution.
In each step, we test our equations at the time $t$ by suitable test functions evaluated at the same time~$t$.
However, we do not write the symbol $t$ for simplicity, and it is understood that the equalities we obtain hold \aet.
For the reader's convenience, we recall the definition of the energy~$\calE$ (see~\eqref{Ienergy}) \an{related to the system}
\begin{align}
  \calE(\phi,a,n,\sig)
  & = \iO a (\ln a - 1)
  - \chiphi \iO n \phi
  - \chia \iO a \sig
  \non
  \\
  & {}\quad + \frac 12 \iO \bigl(
    |\nabla\phi|^2
    + |\nabla n|^2
    + |\nabla\sig|^2
  \bigr)
  + \iO F(\phi)
  \label{energy}
\end{align}
and notice that its time derivative is given~by
\begin{align}
  & \ddt \, \calE(\phi,a,n,\sig)
  \non
  \\
  & = \iO \an{\dt a \, \ln a}
  - \chiphi \iO \an{ \dt n \, \phi}
  - \chiphi \iO n  \, \dt\phi
  - \chia \iO \an{ \dt a\, \sig{}}
  \an{{}  - \chia \iO a \, \dt\sig
}\non
  \\
  & \quad {}
  + \frac 12 \, \ddt \iO \bigl(
    |\nabla\phi|^2
    + |\nabla n|^2
    + |\nabla\sig|^2
  \bigr) 
  + \ddt \iO F(\phi) \,.
  \label{dtenergy}
\end{align}

\step
First a priori estimate

We test \eqref{prima} by
$\mu$ and $-\chiphi n$ \an{to} obtain~\pier{the equalities}
\begin{align}
  & \iO \dt\phi \, \mu
  + \iO |\nabla\mu|^2
  - \chiphi \iO \nabla n \cdot \nabla\mu
  = \iO \bigl( -m\phi + \h(\phi) \bigr)\last{ \mu }\pier{,}
  \non
  \\
  & - \chiphi \iO \dt\phi \, n
  - \chiphi \iO \nabla\mu \cdot \nabla n
  + \chiphi^2 \iO |\nabla n|^2
  = \chiphi \, m \iO \phi n
  - \chiphi \iO \h(\phi) \, n \,.
  \non
\end{align} 
Now, \pier{we recall the definition of $R$ in \eqref{HPinterior} and set}
\Beq
  M := \frac 1{\deltaz} \, (m R + 1) \,,
  \label{defM}
\Eeq
where $\deltaz$ is the same as in \eqref{trickMZ}.
We notice at once that the value of $M$ depends only on the structure of the original system,
so that it can be absorbed in the notation $c$ for the generic constants in performing estimates.
Thus, we \juerg{keep $M$ explicitly} only when it is needed.
\an{Then}, we test \eqref{seconda} by $\dt\phi$, $M(\phi-\phibar)$, and $m\phi-\h(\phi)$ \an{to infer that}
\begin{align}
  & \frac 12 \, \ddt \iO |\nabla\phi|^2
  + \ddt \iO F(\phi)
  = \iO \mu \, \dt\phi \, \pier{,}
  \non
  \\
  & M \iO \an{|\nabla \phi|^2}
  + M \iO F'(\phi) (\phi-\phibar)
  = M \iO \mu (\phi-\phibar) \, \pier{,}
  \non
  \\
  & m \iO |\nabla\phi|^2 
  - \iO \h'(\phi) |\nabla\phi|^2
  + m \iO F'(\phi) \, \phi
  - \iO F'(\phi) \, \h(\phi)
  = \iO \mu \bigl( m\phi - \h(\phi) \bigr) \,.
  \non
\end{align}
Next, we test \eqref{terza} by $\ln a-\chia\sig$.
By noting \last{the identity} $\nabla a-\chia a\nabla\sig=a\nabla(\ln a-\chia\sig)$, we have~that
\Beq
  \iO \dt a \, \ln a
  - \chia \iO \dt a \, \sig
  + \iO a |\nabla(\ln a-\chia\sig)|^2
  = \iO (a-a^2+u) (\ln a-\chia\sig) \,.
  \non
\Eeq
Finally, we test \eqref{quarta} and \eqref{quinta} by $\dt n$ and $\dt\sig$, respectively\an{, leading to}
\begin{align}
  & \iO |\dt n|^2
  + \frac 12 \, \ddt \iO |\nabla n|^2
  -  \chiphi \iO \phi \dt n
  = \iO \Sn \dt n \, \pier{,}
  \non
  \\
  & \iO |\dt\sig|^2
  + \frac 12 \, \ddt \iO |\nabla\sig|^2
  -  \chia \iO a \dt\sig
  = \iO \bigl( (1-\sig) -a\sig \bigr) \dt\sig \,.
  \non
\end{align}
At this point, we add all the above equalities to each other
and \last{to} \juerg{the sides of the resulting identity} the equal terms $\ddt\iO|n|^2$ and $2\iO n\dt n$, respectively.
Notice that four terms cancel \juerg{out} 
and that nine of the contributions on the \lhs\ yield those of the time derivative \eqref{dtenergy} of the energy.
By also rearranging \an{terms}, we conclude~that
\begin{align}
  & \ddt \, \calE(\phi,a,n,\sig)
  + \iO |\nabla\mu|^2
  + \chiphi^2 \iO |\nabla n|^2
  \non
  \\
  & \quad {}
  + M \iO \an{|\nabla\phi|^2}
  + M \iO \Funo'(\phi) (\phi-\phibar)
  + m \iO |\nabla\phi|^2 
  \non
  \\
  & \quad {}
  + m \iO \Funo'(\phi) \, (\phi-\rz)
  + \iO a |\nabla(\ln a-\chia\sig)|^2
  \non
  \\
  & \quad {}
  + \iO |\dt n|^2
  + \iO |\dt\sig|^2
  + \iO a^2 (\ln a-\chia\sig)
  + \ddt \iO |n|^2
  \non
  \\
  \separa
  & = 2 \chiphi \iO \nabla\mu \cdot \nabla n
  - \chiphi \iO \bigl( \h(\phi) - m\phi \bigr) n
  \non
  \\
  & \quad {}
  - M \iO \Fdue'(\phi) (\phi-\phibar)
  + M \iO \mu (\phi-\phibar) 
  \non
  \\
  & \quad {}
  + \iO \h'(\phi) |\nabla\phi|^2
  + \iO \Funo'(\phi) \bigl( -m\rz + \h(\phi) \bigr)
  + \iO \Fdue'(\phi) \bigl( -m\phi + \h(\phi) \bigr)
  \non
  \\
  & \quad {}
  + \iO (a+u) (\ln a-\chia\sig)
  + \iO \Sn \dt n
  \non
  \\
  & \quad {}
  + \iO \bigl( (1-\sig) -a\sig \bigr) \dt\sig 
  + 2 \iO n \dt n \,.
  \label{per1stima}
\end{align}

We consider the terms on the \lhs\ of \eqref{per1stima} that involve~$\Funo'$.
The second one is nonnegative, since \juerg{$\Funo'$} is monotone and vanishes at~$\rz$.
As for the other, we recall \eqref{meanok} and apply \eqref{trickMZ} to obtain~that
\Beq
  M \iO \Funo'(\phi) (\phi-\phibar)
  \geq M \deltaz \iO |\Funo'(\phi)| 
  - c \,.
  \label{daMZ}
\Eeq
Since $\chia\in(0,1)$ and \eqref{stimasig} holds, we have~that
\Beq
  \iO a^2 (\ln a-\chia\sig)
  \geq \iO a^2 (\ln a-1)\,,
  \non
\Eeq
and we notice that the last integrand is bounded from below.

Let us come to the \rhs\ of \eqref{per1stima}, where just some terms need an accurate treatment.
Since $\phi-\phibar$ has zero mean value, 
by also applying the \Poinc\ inequality and recalling \eqref{meanok}, we derive~that
\begin{align}
  & M \iO \mu (\phi-\phibar)
  = M \iO (\mu-\mubar) (\phi-\phibar) 
  \non
  \\
  & \leq \frac 14 \iO |\nabla\mu|^2
  + c \iO |\phi-\phibar|^2
  \leq \frac 14 \iO |\nabla\mu|^2
  + c \iO |\phi|^2  
  + c \,.
  \non
\end{align}
As for the term involving $\Funo'$, we \pier{recall \eqref{HPinterior} and observe}~that
\Beq
  \iO \Funo'(\phi) \bigl( -m\rz + \h(\phi) \bigr)
  \leq \pier{{}\sup_{r\in \erre} |\h (r) - m\rz | \iO |\Funo'(\phi)| = 
   m R \iO |\Funo'(\phi)| {}}\,.
  \non
\Eeq
Thus, due to the choice \eqref{defM} of $M$, this term \an{can be absorbed on} the \lhs.
By \eqref{daMZ}, we have indeed
\begin{align}
  & M \iO \Funo'(\phi) (\phi-\phibar)
  - \iO \Funo'(\phi) \bigl( -m\rz + \h(\phi) \bigr)
  \non
  \\
  & \geq \pier{( M \, \deltaz - \pier{mR} )} \iO |\Funo'(\phi)|  \pier{{}- c{}}
  = \iO |\Funo'(\phi)| \pier{{}- c{}} \,.
  \non
\end{align}
The integral involving $a$ and $u$ is treated by recalling that
\an{$a$, $\sig$,} and $u$ are nonnegative and that $u$ is bounded by~$\umax$.
\an{Namely, we} have~that
\begin{align}
  & \iO (a+u) (\ln a-\chia\sig) 
  \leq \iO (a+u) \ln a
  \leq \frac 14 \iO a^2 (\ln a - 1) + c \,.
  \non
\end{align}
Finally, \an{we observe that}
\begin{align}
  & \iO \bigl( (1-\sig) - a\sig \bigr) \dt\sig 
  \leq \iO (1 + a) |\dt\sig|
  \non
  \\
  & \leq \frac 12 \iO |\dt\sig|^2
  + \iO |a|^2 + c
  \leq \frac 12 \iO |\dt\sig|^2
  + \frac 14 \iO a^2 (\ln a - 1) + c \,.
  \non
\end{align}
By recalling the definition of $\Sn$ given in \eqref{quarta}, 
that both $\Fdue'$ and $\h$ are \Lip\ continuous, and that $\h$ is even bounded,
the other terms on the \rhs\ of \eqref{per1stima}
can easily be treated using Young's inequality.
Hence, collecting all the above estimates and \eqref{per1stima} itself,
and ignoring some nonnegative terms on the \lhs, we conclude~that
\begin{align}
  & \ddt \, \calE(\phi,a,n,\sig)
  + \frac 12 \iO |\nabla\mu|^2
  +  \iO |\Funo'(\phi)|
  + \iO a |\nabla(\ln a-\chia\sig)|^2
  \non
  \\
  & \quad {}
  + \frac 12 \iO |\dt n|^2
  + \frac 12 \iO |\dt\sig|^2
  + \frac 12 \iO a^2 (\ln a-1)
  + \frac 12 \, \ddt \iO |n|^2
  \non
  \\
  & \leq c \iO |\nabla\phi|^2
  + c \iO |\phi|^2 
  + c \iO |\nabla n|^2
  + c \iO |n|^2
  + c \,.
  \non
\end{align}

\an{At this point, we integrate the resulting inequality over $(0,t)$, where $t \in (0,T)$ is arbitrary. The left-hand side contains two terms with no prescribed sign, specifically those of $\calE$ involving the products $n\phi$ and $a\sig$. Consequently, we cannot directly apply the Gronwall lemma.}
\juerg{We therefore move them to} the \rhs\ and estimate them. 
To this end, recall that \pier{$\chiphi, \, \chi_a \in (0,1)$} and that $0\leq\sig\leq1$.
Moreover, we observe~that 
\Beq
  r^2 \leq \frac 12 \, F(r) + c
  \quad \hbox{for every $r\in D(F)$}.
  \label{daHPF}
\Eeq
\an{T}his \an{readily} follows from \eqref{HPregpot} in the case of regular potentials,
and it is trivially satisfied in the case of potentials satisfying~\eqref{HPlogpot}.
Hence, we have~that
\begin{align}
  & \chiphi \iO n(t) \phi(t)
  + \chia \iO a(t) \sig(t)
  \leq \frac 14 \iO |n(t)|^2 
  + \iO |\phi(t)|^2
  + \iO a(t)
  \non
  \\
  & \leq \frac 14 \iO |n(t)|^2 
  + \frac 12 \iO F(\phi(t))
  + \frac 12 \iO a(t)(\ln a(t) - 1)
  + c \,,
  \non
\end{align}
and this \an{can be absorbed on} the \lhs.
Moreover, there is one more term on the \rhs\ to be treated,
namely, the term \an{arising} from the time integration of $\iO|\phi|^2$.
\an{Th}is \an{latter} can be estimated by using \eqref{daHPF} once more.
\juerg{Now recall} that our assumptions on $F$ imply that $F$ is bounded from below\juerg{. We therefore}
can apply Gronwall's lemma\an{, \pier{owing also to} the assumptions on the initial data,} and conclude~that
\begin{align}
  & \norma\phi_{\L\infty V}
  + \norma{F(\phi)}_{\L\infty\Luno}
  + \norma{\nabla\mu}_{\L2H}
  \non
  \\
  & \quad {}
  + \norma{a(\ln a-1)}_{\L\infty\Luno}
  + \norma{a^2(\ln a-1)}_{\LQ1}
  + \norma{a^{1/2}\nabla(\ln a - \chia\sig)}_{\LQ2}
  \non
  \\
  & \quad {}
  + \norma n_{\H1H\cap\L\infty V}
  + \norma\sig_{\H1H\cap\L\infty V}
  \leq c \,.
  \label{1stima}
\end{align}
\pier{The bound in \eqref{1stima}} also implies that
\Beq
  \norma a_{\L\infty\Luno\cap\LQ2} \leq c \,.
  \label{stimaa}
\Eeq

\step
Consequences

\pier{At this point, it is rather straightforward to infer that
\Beq
  \norma{\dt\phi}_{\L2\Vp} + \norma n_{\L2W} + \norma\sig_{\L2W} \leq c \,.
  \label{nsigdtphi}
\Eeq
Indeed,} for the first of these estimates, one tests \eqref{prima} by any $v\in\L2V$, 
integrates over~$(0,T)$, and accounts for \eqref{1stima}, to obtain~that
\Beq
  \ioT \< \dt\phi(t) , v(t) > \, dt
  \leq c \, \norma v_{\L2V} \,.
  \non
\Eeq
\pier{By rewriting \eqref{quarta} and \eqref{quinta} as partial differential equations 
\begin{align}
  &  \dt n 
-\Delta  n 
  = \chiphi \phi 
  + \cphi \phi + \cn n + \csig \sig + c_0 \quad\aeQ \,,
  \label{quarta-pde}
  \\
  & \dt\sig - \Delta \sig = (1-\sig) + a (\chia-\sig) \quad\aeQ \,,
  \label{quinta-pde}
\end{align}
and comparing the terms, \jurgen{we infer from \eqref{1stima}, \eqref{stimaa} and \eqref{stimasig}} that 
\Beq
   \norma{\Delta n}_{\L2H} + \norma{\Delta \sig}_{\L2H} \leq c \,,
  \non
\Eeq
whence we obtain \eqref{nsigdtphi} with the help of the} elliptic regularity \an{theory}.

\step
Second a priori estimate

We test \eqref{seconda} by $\phi-\phibar$
and owe to \eqref{daMZ} divided by~$M$.
By also accounting for the \Poinc\ and Young inequalities and \eqref{meanok}, we obtain~that
\begin{align}
  & \deltaz \iO |\Funo'(\phi)|
  \leq \iO \Funo'(\phi) (\phi-\phibar) + c
  \leq \iO |\nabla\phi|^2
  + \iO \Funo'(\phi) (\phi-\phibar) + c
  \non
  \\
  &\quad{} = - \iO \Fdue'(\phi) (\phi-\phibar)
  + \iO \mu (\phi-\phibar)
  + c
  \non
  \\
  &\quad{}= - \iO \Fdue'(\phi) (\phi-\phibar)
  + \iO (\mu-\mubar) (\phi-\phibar)
  + c
  \non
  \\
  &\quad{} \leq c \iO |\phi|^2
  + c \, \norma{\juerg{\nabla\mu}} \, \norma{\phi-\phibar}
  + c \,.
  \label{per2stima}
\end{align}
Now we square, integrate over~$(0,T)$, and apply \eqref{1stima}.
This yields~that
\Beq
  \norma{\Funo'(\phi)}_{\L2\Luno} \leq c\,,
  \quad \hbox{whence} \quad
  \norma{\overline{\Funo'(\phi)}}_{L^2(0,T)} \leq c \,.
  \non
\Eeq
Therefore, by testing \eqref{seconda} by $1/|\Omega|$ and comparing,
we infer~that
\Beq
  \norma\mubar_{L^2(0,T)} \leq c \,.
  \non
\Eeq
By combining with \eqref{1stima}, \pier{and using once more inequality~\eqref{poincare},} we conclude~that
\Beq
  \norma\mu_{\L2V} \leq c \,.
  \label{2stima}
\Eeq

\step 
Consequence

Now, we consider \eqref{seconda}.
By splitting $F$ as $F=\Funo+\Fdue$, moving the term involving $\Fdue'$ to the \rhs\
and applying a usual argument based on the monotonicity of~$\Funo'$ \pier{(i.e., one can test \eqref{seconda}
by $\Funo' (\phi)$), we deduce~that both $\Funo'(\phi)$ and $\Delta\phi$} are estimated in $H$ by the $H$ norm of the \rhs.
Therefore, \jurgen{we conclude from elliptic regularity} that
\Beq
  \norma\phi_{\L2W} + \norma{\Funo'(\phi)}_{\L2H} \leq c \,.
  \label{phiW}
\Eeq

\step
Third a priori estimate

\pier{We take $v=a$ in \eqref{terza} and, using the positivity of $a$ and Young's inequality, we have that 
\begin{align}
  & \frac 12 \, \ddt \norma{a}^2 
  +  \iO |\nabla a|^2 
  + \iO |a|^3 
  \non\\
  &  
  \leq \chia\iO a \nabla\sig \cdot \nabla a + \frac 3 2 \norma{a}^2 + \frac 1 2  \last{\norma{\umax}_\infty^2} |\Omega|. 
\label{pier3} 
\end {align} 
Now, in order to deal with the first integral on the \rhs, we invoke H\"older's inequality and the \Lady\ inequalities in \eqref{lady}, and find out that 
\begin{align}
  & \chia\iO a \nabla\sig \cdot \nabla a \leq  \norma{a}_4 \norma{\nabla \sigma}_4
   \norma{\nabla a}_2
   \non\\
   &\quad{}
\leq c \norma{a}^{1/2} \norma{a}_V^{1/2}  \norma{\sigma}_V^{1/2}  (\norma \sigma + \norma{\Delta \sigma})^{1/2} \norma{a}_V
  \non\\
   &\quad{}
\leq c \norma{a}^{1/2}  \norma{\sigma}_{\L\infty V}^{1/2}  \norma{\sigma}_{W}^{1/2} \norma{a}_V^{3/2}
  \non\\
   &\quad{}
\leq \frac 1 2 \norma{a}_V^{2} + c \norma{\sigma}_{W}^{2} \norma{a}^2 ,
\label{pier4} 
\end {align}
where \eqref{1stima} and the Young inequality~\eqref{young}, with exponents $4/3$ and $4$, have been exploited. Note that the \last{function} $t\mapsto \norma{\sigma (t)}_{W}^{2}$ is known to be bounded in $L^1(0,T)$ by \eqref{nsigdtphi}. Then, combining \eqref{pier3} and \eqref{pier4}, we can integrate the resultant over $(0,t)$ with the help of the initial condition for $a$, see \eqref{cauchy} and \eqref{HPaz}. 
Next, we apply the Gronwall lemma and deduce that}
\Beq
  \pier{\norma a_{\L\infty H\cap\L2V\cap\LQ3} \leq c\,,
  \quad \hbox{whence also (cf.~\eqref{disug})} \quad
  \norma a_{\LQ4} \leq c \,.}
  \label{3stima}
\Eeq

\step 
\pier{Regularity for two variables}

We set $g:=(1-\sig)+a(\chia-\sig)$ for a while.
Then, \pier{we have that} $g\in\LQ4$ \pier{by \eqref{stimasig} and~\eqref{3stima}.}
Moreover, since
\Beq
  \nabla g = \juerg{-}\nabla\sig + \nabla a \last{(\chia-\sig) }- a \nabla\sig
  \non
\Eeq
and $\nabla\sig\in\LQ4$ by \eqref{1stima}, \eqref{nsigdtphi}, and~\eqref{embeddings},
we see that $\nabla g\in\LQ2$, so that $g$ belongs to $\L2V$
\pier{as well. We also recall the property $\sigma_0 \in W$ from~\eqref{HPsigz}.
Therefore, on \jurgen{the} one hand, a comparison in \eqref{quinta-pde} and maximal parabolic regularity 
 (see \cite[Thm.~2.1]{DHP}) yield~that 
\Beq
  \norma\sig_{\W{1,4}{\Lx4}\cap\L4{\Wx{2,4}}}
  \leq c \,.
  \label{maxregpsig}
\Eeq
On the other hand, by the abstract regularity theory for parabolic problems contained, e.g., in \cite{Lions},
we also infer that
\Beq
  \norma\sig_{\H1{V}\cap\C0W\cap\L2{\Hx{3}}}
  \leq c \,.
  \label{maxregHsig}
\Eeq
In view of \eqref{HPsigz} and} \eqref{quarta-pde}, a similar conclusion holds for $n$ since the \pier{\rhs\ $(\chi_\phi + \cphi) \phi + \cn n + \csig \sig + c_0 $ is} bounded both in $\LQ4$ and in $\L2V$ due to \pier{\eqref{1stima}.}
Th\last{us}, we have that
\begin{align}
\pier{\norma{n}_{\W{1,4}{\Lx4}\cap\L4{\Wx{2,4}}}+
  \norma{n}_{\H1{V}\cap\C0W\cap\L2{\Hx{3}}}
  \leq c \,.}
  \label{maxregn}
\end{align}

\step
Fourth a priori estimate

\pier{From \eqref{terza} and integration by parts\last{,} it follows that 
\Beq   \iO \dt a \, v
  + \iO \nabla a \cdot \nabla v
  = - \chia \iO  (\nabla a \cdot \nabla \sigma +   a \Delta \sig ) v
  + \iO \bigl( a - a^2 + u \bigr) v  
\label{pier5}
\Eeq
for every $v\in V$, a.e.~in $(0,T)$. Note that the first integral on the \rhs\ makes sense since $\nabla \sigma $ is bounded in $\L4{\Lx\infty}$ by \eqref{maxregpsig} and the Sobolev embedding $\Wx{2,4}\subset 
\Lx\infty$, whereas $a$ and $\Delta \sig$ are bounded in $\LQ4$. We formally take $v = \dt a $ in \eqref{pier5} and, by the \Holder\ and Young inequalities, we easily obtain
\begin{align}
  &  \norma{\dt a}^2 
  + \frac 12 \, \ddt \|\nabla a\|^2 
  \non\\
  &  
  \leq \norma{\nabla a} \norma{\nabla\sig}_\infty  \norma{\dt a} 
  + \norma{a}_4 \norma{\Delta \sig}_4  \norma{\dt a}
   + \norma{a-a^2 +u}_2 \norma{\dt a}
     \non\\
  &\leq \frac 12 \norma{\dt a}^2 + c  \norma{\sig}_{\Wx{2,4}}^2 \norma{\nabla a}^2 
  + c  \norma{a}_{\Lx4}^2 \norma{\sig}_{\Wx{2,4}}^2 + c \bigl(\norma{a}_{\Lx4}^4 +1\bigr).
\label{pier6} 
\end {align} 
Then, in view of \eqref{3stima} and \eqref{maxregpsig}, we are allowed to integrate over $(0,t)$ and apply the Gronwall lemma \last{as $t\mapsto \norma{\sig (t)}_{\Wx{2,4}}^2$ is bounded in $L^2(0,T)$}
to conclude that 
\Beq
  \norma{a}_{\H1H\cap \L\infty V} \leq c \,.
  \label{quartastima}
\Eeq
Then, going back to \eqref{pier5} and emphasizing that now the whole term $\nabla a \cdot \nabla \sigma +   a \Delta \sig$ is bounded in $\L2H$, by comparison and elliptic regularity we infer that}
\Beq
 \pier{ \norma{a}_{\L2 W} \leq c \,.}
  \label{4stima+}
\Eeq
 
%

\step
Fifth a priori estimate

\pier{We proceed formally and take $v = \dt \mu $ in \eqref{prima}. Coincidently, we differentiate \eqref{seconda} 
and test the resulting equation by $\dt \phi$. Then, we sum up, noting that a cancellation occurs, and integrate also by parts over $(0,t)$. We obtain
\begin{align}
  & \frac 12 \iO |\nabla\mu(t)|^2
  + \intQt |\nabla\dt\phi|^2
  + \intQt F_1''(\phi)|\dt\phi|^2
  \non
  \\
  & = \frac 12 \iO |\nabla\mu(0)|^2
+ \chi_\phi\iO \nabla n (t)\cdot \nabla \mu (t)   
- \chi_\phi \iO \nabla n_0 \cdot \nabla \mu (0) 
\non
  \\
  & \quad{}
  -\intQt \nabla \dt n \cdot \nabla \mu   
+ \iO (\h(\phi) - m\phi)  (t) \mu (t)   
- \iO (\h(\phiz) - m\phiz) \mu (0)   
\non
  \\
  & \quad{}
- \intQt (\h' (\phi) - m )  \dt\phi\, \mu   
- \intQt F_2''(\phi) |\dt\phi|^2,
  \label{pier7}
\end{align}
where we \last{employed} the notation $Q_t= \Omega\times (0,t).$ 
We point out that the third term on the \lhs\ is nonnegative due to the monotonicity of $F_1'$. About the value $\mu(0)$, we recover it from \eqref{seconda} and realize from \eqref{HPphiz} that $\mu(0)= \mu_0$ is bounded in $V$. Then, recalling also \eqref{HPsigz} and \eqref{HPh}, it is clear that the first, third and sixth terms on the \rhs\ of \eqref{pier7} are under control. The second one can be easily treated as 
$$ \chi_\phi\iO \nabla n (t)\cdot \nabla \mu (t)  \leq \frac 18 \iO |\nabla\mu(t)|^2
+ c \norma{n}_{\L\infty V}^2 \leq \frac 18 \iO |\nabla\mu(t)|^2 +c $$
by the Young inequality and \eqref{maxregn}. The fourth term on the \rhs\ of \eqref{pier7} is already bounded due to \eqref{maxregn} and \eqref{2stima}. On the other hand, as the functions $\h'$ and $F_2''$ are bounded, by virtue of the Young inequality, the compactness inequality~\eqref{compact}, and \eqref{nsigdtphi}, we deduce that
\begin{align*}
&{}- \intQt (\h' (\phi) - m )  \dt\phi\, \mu   
- \intQt F_2''(\phi) |\dt\phi|^2 \leq c + c \intQt|\dt\phi|^2  
\\
&\quad{}
\leq c + \frac 12  \intQt |\nabla\dt\phi|^2
+ \norma{\dt\phi}_{\L2\Vp}^2 \leq c + \frac 12  \intQt |\nabla\dt\phi|^2. 
\end{align*} 
Now, in \eqref{pier7} it remains to control one term, for which we use the boundedness of $\h$ and the estimate~\eqref{1stima}, along with the \Poinc\ inequality~\eqref{poincare} and, once more, the equation~\eqref{seconda} with $v=1/|\Omega|$. It follows that 
\begin{align*}
 &\iO (\h(\phi) - m\phi)  (t) \mu (t)   \leq c\norma {\mu (t)}  \leq c\norma {\mu (t) -\ov \mu (t)} + c \bigl| \ov \mu (t)\bigr| \\
 &\quad{} \leq c \norma {\nabla \mu (t)} + c \iO |\Funo'(\phi(t))| + c \norma{\phi (t)} + c .
\end{align*}
Now, we recall \eqref{per2stima} and arrive at 
\begin{align*}
&\iO (\h(\phi) - m\phi)  (t) \mu (t)  \\
 &\quad{}  \leq 
 c \norma {\nabla \mu (t)} \bigl( 1+ \norma{\phi(t)-\phibar(t)} \bigr) + 
c \norma{\phi(t)}^2
  + c \leq  \frac 18 \iO |\nabla\mu(t)|^2 +c .
\end{align*}
Then, collecting the above computations in \eqref{pier7} and invoking \eqref{nsigdtphi}
lead to the estimate 
\Beq
 \norma{\nabla\mu}_{\L\infty H} + \norma{\dt\phi}_{\L2V} \leq c \,.
  \label{pier8}
\Eeq
Hence, recalling \eqref{per2stima} again, at this point we can infer that 
\Beq
  \norma{\Funo'(\phi)}_{\L\infty\Luno} \leq c\,,
  \non
\Eeq
whence, testing \eqref{seconda} by $1/|\Omega|$ and comparing, it holds that
$\, 
  \norma\mubar_{L^\infty (0,T)} \leq c
$,
and, consequently, using once more inequality~\eqref{poincare}, we conclude~that
\Beq
  \norma\mu_{\L\infty V} \leq c \,.
  \label{pier9}
\Eeq
Moreover, by \eqref{pier8} and a comparison of terms in the \last{strong form} of
\eqref{prima}, i.e., 
\Beq
\dt\phi  - \Delta (\mu - \chiphi n) 
  = - m \phi + \h(\phi) \quad \hbox{a.e. in } \, Q,
\label{prima-pde} 
\Eeq
we find out that $\Delta (\mu - \chiphi n) $ is bounded in $\L{2}V$, whence, by elliptic regularity, 
$\mu - \chiphi n$ is bounded in $\L{2}{W\cap\Hx3}$ and, consequently, in view of \eqref{maxregn}, 
it holds that
\Beq
 \norma\mu_{\L{2}{W\cap\Hx3}} \leq c \,.
  \label{pier9bis}
\Eeq
Now, we can argue in the same way as for \eqref{phiW} and find as well  that 
\Beq
  \norma\phi_{\L\infty W} + \norma{\Funo'(\phi)}_{\L\infty H} \leq c \,.
  \label{pier10}
\Eeq
The information given by the estimates~\eqref{pier8} and \eqref{pier10} is enough to conclude that $\phi\in C^0(\ov Q) $ so that the values assumed by $\phi$ range in a compact subset of $\erre$, and in the case $D(F)=\erre$  we have already proved the separation property \eqref{separation}.}

\step
Further estimate

\pier{From now on, we restrict ourselves to the case $D(F)=(0,1)$ and argue similarly as in the proof of  \cite[Proposition~2.6]{CGRS} in order to obtain the estimate 
\Beq
  \norma\phi_{\W{1,\infty}H\cap \H1W } +  \norma{\mu}_{\LQ\infty }\leq c \,.
  \label{pier11}
\Eeq
In fact, the technique is as follows. From \eqref{prima}, considered at the initial time $t=0$, and
\eqref{HPphiz}, \eqref{HPsigz} we recover that 
$$ \dt \phi (0) = \Delta (\mu_0 -\chi_\phi n_0) - m \phiz +\h(\phiz) $$
is bounded in $H$. Next, we differentiate both \eqref{prima} and \eqref{seconda} and test by $\dt \phi$ and $-\Delta (\dt \phi)$, respectively. Then, we sum up and integrate by parts and over $(0,t),$ paying attention to the cancellation of two terms.  The integrals on the \rhs\  
$$ 
\chi_\phi\intQt \nabla \dt n \cdot \nabla \dt \phi  
+ \intQt (\h' (\phi) - m) |\dt \phi|^2   
$$
\last{are readily controlled}
thanks to \eqref{maxregn} and \eqref{pier8}.
Then we arrive at some inequality similar to \cite[formula~(5.16) and following]{CGRS}. From this point we can proceed along the same line as in \cite[pp.~2160-2162]{CGRS}, by exploiting the inequality in \eqref{HPlogpot}, which gives a control for $F_1''$, along with the Trudinger inequality (see, e.g., \cite{Moser})
\Beq
  \iO e^{|v|}
  \leq c_\Omega \, e^{c_\Omega\,\normaV v^2}
  \quad \hbox{for every $v\in V$}\,,
  \non
\Eeq
which holds in \last{two dimensions}. 
Let us omit the full proof here: this proof permits to arrive at the estimate 
\Beq
  \norma{\dt\phi}_{\L{\infty}H} +  \norma{\Delta (\dt \phi)}_{\L2 H}\leq c ,
  \non
\Eeq
from which by the elliptic regularity theory and \eqref{pier8}, \eqref{pier10} we achieve the estimate \eqref{pier11} for the term involving $\phi.$ \last{Next, a 
comparison argument 
in \eqref{prima-pde} reveals}
that $\Delta (\mu - \chiphi n) $ is bounded in $\L{\infty}H$, whence by elliptic regularity 
and \eqref{maxregn}, $\mu$ is shown to be bounded in $\L{\infty}W \subset \LQ\infty $, which completes the proof of \eqref{pier11}.}

\step
Separation property

Up to now, we have \pier{completely proved the stability estimate \eqref{stability}, 
which} holds true with a constant $K_1$ satisfying the properties given in the statement,
since the constants $c$ we have introduced \pier{in the various steps enjoy these properties. 
We still have to check the separation property~\eqref{separation} in the case $D(F)=(0,1)$.
To this aim, 
we start from the bound for $\mu$ in~$\LQ\infty$ in \eqref{pier11}.}
Then, a Moser type procedure in \eqref{seconda} provides a bound for~$\Funo'(\phi)$.
Let us sketch this argument by proceeding formally \pier{and}
\an{acknowledging that a truncation argument would suffice to obtain a rigorous proof}.
To simplify notation,  we set $\psi:=\Funo'(\phi)$ and $g:=\mu-\Fdue'(\phi)$.
The argument just deals with the elliptic equation at the time~$t$,
which however is not written for simplicity.
We take any $p>2$ and test \eqref{seconda} by $|\psi|^{p-2}\psi$ to obtain~that
\Beq
  (p-1) \iO |\psi|^{p-2} |\nabla\phi|^2
  + \iO |\psi|^p
  = \iO g |\psi|^{p-2}\psi \,.
  \non
\Eeq
\pier{By the Young inequality~\eqref{young}} we deduce that
\Beq
  \norma\psi_p^p
  \leq \iO |g| \, |\psi|^{p-1}
  \leq \norma g_p \, \norma{|\psi|^{p-1}}_{p'}
  = \norma g_p \, \norma\psi_p^{p/p'} 
  \leq \frac 1p \, \norma g_p
  + \frac 1{p'} \, \norma\psi_p^p \,.
  \non
\Eeq
By rearranging, we infer that $\norma\psi_p\leq\norma g_p$.
Then, letting $p$ tend to infinity, we conclude that $\norma\psi_\infty\leq\norma g_\infty$\pier{,
which entails 
$\,\norma{\Funo'(\phi(t))}_\infty
  \leq \norma{\mu(t)-\Fdue'(\phi(t))}_\infty \, $
\aat, whence} 
$$ \norma{\Funo'(\phi)}_{\LQ\infty}
  \leq \norma{\mu-\Fdue'(\phi)}_{\LQ\infty}
  \leq c \,.
$$
By accounting for assumption \eqref{HPlogpot}, we deduce that \eqref{separation} holds true
with some values $r_\pm$ as in the statement.
\qed

\section{Uniqueness and continuous dependence}
\label{UNIQUENESS}
\setcounter{equation}{0}

In this section, we prove the uniqueness part of Theorem~\ref{Wellposedness} and the continuous dependence estimate~\eqref{contdep}.
\an{We} just prove the latter for arbitrary solutions corresponding to the control variables,
so that uniqueness follows as a consequence of the case of the same control.
\an{In this direction}, we fix \an{$u_i\in \Uad$, $i=1,2$,}
and any two corresponding solutions $(\phi_i,\mu_i,a_i,n_i,\sig_i)$ with the regularity 
and the properties stated in Theorem~\ref{Wellposedness}, in particular, the separation property~\eqref{separation}.
We set for convenience
\begin{align*}
  &u:= u_1 - u_2 \,, \quad
  \phi = \phi_1 - \phi_2 \,, \quad
  \pier{\mu = \mu_1 - \mu_2 \,,} \\
  &
  a = a_1 - a_2 \,, \quad
  n:= n_1 - n_2\,, \quad 
  \sig:= \sig_1 - \sig_2 \,.
  \non
\end{align*}
Then, we write all the equations \accorpa{prima}{quinta} for both solutions  and take the differences \an{to} obtain~that
\begin{align}
  & \iO \dt\phi \, v
  + \iO \nabla\mu \cdot \nabla v
  - \chiphi \iO \nabla n \cdot \nabla v
  = - m \iO \phi v
  + \iO \bigl( \h(\phi_1) - \h(\phi_2) \bigr) \, v\,,
  \label{dprima}
  \\
  & \iO \nabla\phi \cdot \nabla v
  + 
  \iO \bigl( F'(\phi_1) - F'(\phi_2) \bigr) \, v
  = \iO \mu v \,,
  \label{dseconda}
  \\
  & \pier{\iO \dt a \, v}
  + \iO \nabla a \cdot \nabla v
  - \chia \iO \bigl( a \nabla\sig_1 + a_2 \nabla\sig \bigr) \cdot \nabla v
  = \iO \bigl( a - (a_1+a_2)a + u \bigr) v \,,
  \label{dterza}
  \\
  & \iO \dt n \, v
  + \iO \nabla n \cdot \nabla v
  - \chiphi \iO \phi v
  = \iO \bigl( \cphi\phi + \cn n + \csig\sig \bigr) \, v \,,
  \label{dquarta}
  \\
  & \iO \dt\sig \, v
  + \iO \nabla\sig \cdot \nabla v
  = \iO \bigl(  -\sig + \chia a - a\sig_1 - a_2\sig \bigr) v\,,
  \label{dquinta}
\end{align}
\juerg{for} every $v\in V$ and \aet.

\Brem
\label{Remcontdep}
In our computations, we allow the values of the generic constants $c$ to also depend on the solutions we are considering.
However, at the end of the proof, since uniqueness follows as a consequence as already said before,
one realizes that the solutions taken into account are exactly those provided by the \juerg{already proved} existence part of 
Theorem~\ref{Wellposedness}.
This implies that the norms of the solutions considered in the present proof are bounded by the constant $K_1$ of the stability estimate~\eqref{stability},
and thus they depend only on $\Omega$, $T$, the structure of the system, the initial data, and~$\umax$.
\Erem

\step
First estimate

We test the above equations by $\phi$, $-\Delta\phi$, $a$, $\dt n$, and~$\dt\sig-\Delta\sig$, respectively.
More precisely, we test \eqref{dprima} and \eqref{dterza} as said,
while we write \eqref{dseconda}, \eqref{dquarta} and \eqref{dquinta} in their strong form,
multiply them by $-\Delta\phi$, $\dt n$, and~$\dt\sig-\Delta\sig$,
and integrate over~$\Omega$.
Notice that the regularity \Regsoluz\ for both solutions allows this procedure.
After some rearrangement, we obtain~that
\begin{align} 
  & \frac 12 \, \ddt \iO |\phi|^2
  + m \iO |\phi|^2
  = - \iO \nabla\mu \cdot \nabla\phi
  + \chiphi \iO \nabla n \cdot \nabla\phi
  + \iO \bigl( \h(\phi_1) - \h(\phi_2) \bigr) \, \phi \,,
  \non
  \\
  & \iO |\Delta\phi|^2
  = - \iO \bigl( F'(\phi_1) - F'(\phi_2) \bigr) \, (-\Delta\phi)
  + \iO \nabla\mu \cdot \nabla\phi \,,
  \non
  \\
  & \frac 12 \, \ddt \iO |a|^2 
  + \iO |\nabla a|^2
  = \chia \iO \bigl( a \nabla\sig_1 + a_2 \nabla\sig \bigr) \cdot \nabla a
  + \iO \bigl( a^2 - (a_1+a_2)a^2 + ua \bigr)\,,
  \non
  \\
  & \iO |\dt n|^2
  + \frac 12 \, \ddt \iO |\nabla n|^2
  = \chiphi \iO \phi \dt n
  + \iO \bigl( \cphi\phi + \cn n + \csig\sig \bigr) \, \dt n\,,
  \non
  \\
  & \iO |\dt\sig|^2
  + \ddt \iO |\nabla\sig|^2
  + \iO |\Delta\sig|^2
  = \iO \bigl( -\sig + \chia a - a\sig_1 - a_2\sig \bigr) (\dt\sig - \Delta\sig)\,.
  \non
\end{align}
At this point, we take the sum of these identities.
Moreover, we add 
$$ \pier{\frac12 \ddt\iO|n|^2 \ \hbox{ and } \ \ddt\iO|\sig|^2} $$ 
to the \lhs\ of the \juerg{resulting equality} 
and the same terms, written in the form $\iO n\/ \dt n$ and $2\iO\sig\/\dt\sig$, to its \rhs.
The new \lhs\ \juerg{then becomes}
\begin{align} 
  & \frac 12 \, \ddt \iO |\phi|^2
  + m \iO |\phi|^2
  + \iO |\Delta\phi|^2
  + \frac 12 \, \ddt \iO |a|^2 
  + \iO |\nabla a|^2
  \non
  \\
  & \quad {}
  + \iO |\dt n|^2
  + \frac 12 \, \ddt \iO \bigl(\an{|n|^2+|\nabla n|^2}\bigr)
  \non
  \\
  & \quad {}
  + \iO |\dt\sig|^2
  + \ddt \iO \bigl(\an{|\sig|^2 + |\nabla\sig|^2}\bigr)
  + \iO |\Delta\sig|^2 \,,
  \non
\end{align}
and we have to estimate the terms on the corresponding \rhs.
However, two of them cancel each other, and most of the others can be simply dealt with by means of Young's inequality,
possibly on account of the \Lip\ continuity of both $\h$ and $F'$ in the interval $[r_-,r_+]$,
considering that we want to apply (after time integration) the Gronwall lemma.
Hence, we only estimate the terms that need some treatment.
On account of \eqref{elliptic}, we have~that
\Beq
  \chiphi \iO \nabla n \cdot \nabla\phi
  \leq \norma{\nabla\phi} \, \norma{\nabla n}
  \leq \frac 14 \iO \bigl( |\phi|^2 + |\juerg{\Delta}\phi|^2 \bigr)
  + c \iO |\nabla n|^2 \,.
  \non
\Eeq
As for the terms involving the product of three factors, we first \pier{employ \eqref{regn} to find} that
\begin{align}
  & \chia \iO a \nabla\sig_1 \cdot \nabla a
  \leq \norma a_4 \norma{\nabla\sig_1}_4 \, \norma{\nabla a}_2
  \non
  \\
  & \leq \frac 18 \iO |\nabla a|^2
  + c \, \juerg{\|\sig_1\|_{\L\infty W}^2} \, \norma a \, (\norma a + \norma{\nabla a})
  \leq \frac 14 \iO |\nabla a|^2
  + c \, \juerg{\norma a^2}  \,.
  \non
\end{align}
\juerg{Moreover, using the Sobolev inequality in \eqref{disug} and the second \Lady\ inequality \pier{in} \eqref{lady} after Young's inequality,} \pier{we have that}
\begin{align}
  & \chia \iO a_2 \nabla\sig \cdot \nabla a
  \leq \norma{a_2}_4 \norma{\nabla\sig}_4 \, \norma{\nabla a}_2
  \non
  \\
  & \leq \frac 14 \iO |\nabla a|^2 
  + c \, \norma{a_2}_4^2 \, \normaV\sig (\norma\sig + \norma{\Delta\sig})
  \non
  \\
  & \leq \frac 14 \iO |\nabla a|^2 
  + c \, \norma{a_2}_4^2 \, \normaV\sig^2
  + \frac 14 \iO |\Delta\sig|^2
  + c \, \norma{a_2}_4^4 \, \normaV\sig^2\,,
  \non
\end{align}
and we observe that the function $t\mapsto\norma{a_2(t)}_4$ belongs to $L^4(0,T)$
(cf.~\eqref{rega} \pier{and \eqref{disug}}). \juerg{In addition}, it is clear that
\Beq
  - \iO (a_1+a_2) a^2
  \leq \norma{a_1+a_2}_4 \, \norma a_4 \, \norma a_2
  \leq \frac 14 \, \normaV a^2
  + c \, \norma{a_1+a_2}_4^2 \, \norma a_2^2 \,,
  \non
\Eeq
\juerg{where the function $t\mapsto\norma{(a_1+a_2)(t)}_4$ belongs to $L^4(0,T)$.}
The terms involving $a\sig_1$ and $a_2\sig$ are treated in a similar way.
Finally, \pier{from Young's inequality it follows that}
\Beq
  \iO ua
  \leq \frac12 \iO |a|^2 + \frac12 \iO |u|^2 \,.
  \non
\Eeq 
Therefore, collecting all these inequalities and those we have omitted,
we estimate the \rhs\ we are considering in a form that is suitable for the application of the Gronwall lemma.
\an{After time integration, using the assumptions on the initial data as well as elliptic regularity theory, we} conclude~that
\begin{align}
  & \norma\phi_{\L\infty H\cap\L2W}
  + \norma a_{\L\infty H\cap\L2V}
  + \norma n_{\H1H\cap\L\infty V}
  \non
  \\
  & \quad {}
  + \norma\sig_{\H1H\cap\L\infty V\cap\L2W}
  \leq c \, \norma u_{\L2H} \,.
  \label{1stimad}
\end{align}
Notice that this implies that
\Beq
  \norma a_{\LQ4} + \norma{\nabla\sig}_{\LQ4}
  \leq c \, \norma u_{\L2H}\,,
  \label{1stimadbis}
\Eeq
thanks to the third inequality \an{in} \eqref{disug}.

\step
Consequence

A comparison \an{argument} in \an{equations} \eqref{dquarta} and \eqref{dseconda}\an{, along with elliptic regularity theory,} thus yields that
\Beq
  \norma n_{\L2W} \leq c \, \norma u_{\L2H}
  \aand
  \norma\mu_{\L2H} \leq c \, \norma u_{\L2H} \,.
  \label{da1stimad}
\Eeq

\step
Second estimate

We test \eqref{dprima} by $\mu$ and \eqref{dseconda} by $\dt\phi$
and add the resulting equalities.
Noting a cancellation, we have~that
\begin{align}
  & \iO |\nabla\mu|^2
  + \frac 12 \, \ddt \iO |\nabla\phi|^2
  = \chiphi \iO \nabla n \cdot \nabla\mu
  - m\iO \phi \mu
  \non
  \\
  & \quad {}
  + \iO \bigl( \h(\phi_1) - \h(\phi_2) \bigr) \mu
  - \iO \bigl( F'(\phi_1) - F'(\phi_2) \bigr) \dt\phi \,,
  \label{pier12}
\end{align}
and just the last term needs some attention \an{as the others can be easily controlled using Young's inequality}.
To deal with it, we recall the \Lip\ continuity of $F$ in the interval $[r_-,r_+]$
and test \eqref{dprima} by $-(F'(\phi_1)-F'(\phi_2))$.
We obtain~that
\begin{align}
  & - \iO \bigl( F'(\phi_1) - F'(\phi_2) \bigr) \dt\phi
  \non
  \\
  & = \iO \nabla\mu \cdot \nabla\bigl( F'(\phi_1) - F'(\phi_2) \bigr)
  \juerg{\,+\,} \chiphi \iO \Delta n \, \bigl( F'(\phi_1) - F'(\phi_2) \bigr) 
  \non
	\\
	&\juerg{\quad\, +\iO \bigl(m\phi-(\h(\phi_1)-\h(\phi_2))\bigr)\,\bigl(F'(\phi_1) - F'(\phi_2) \bigr)}\,.
	\label{pier13}
\end{align}
\juerg{In view of \eqref{1stimad} and \eqref{da1stimad}, and using Young's inequality and the Lipschitz continuity
of both $\h$ and $F'$ on $[r_-,r_+]$, we easily conclude that the sum of the second and third terms on the \rhs\ is bounded by
\,\,$c\|u\|^2_{L^2(0,T;H)}$.}{} 
Regarding the first term on the \rhs, 
we recall the regularity of $\phi_1$ and $\phi_2$ 
given by Theorem~\ref{Wellposedness} and \eqref{nablaphibdd}, \last{to infer}~that
\begin{align}
  & \iO \nabla\mu \cdot \nabla\bigl( F'(\phi_1) - F'(\phi_2) \bigr)
  = \iO \nabla\mu \cdot
    \an{\Bigl(}
      \bigl( F''(\phi_1) - F''(\phi_2) \bigr) \nabla\phi_1
      + F''(\phi_2) \nabla\phi
    \an{\Bigr)}
  \non
  \\
  & \leq c \, \norma{\nabla\mu}_2 \,
    \bigl( 
      c \, \norma\phi_2 \, \norma{\nabla\phi_1}_\infty
      + \norma{F''(\phi_2)}_\infty \, \norma{\nabla\phi}_2
    \bigr) 
  \leq c \, \norma{\nabla\mu} \, (\norma\phi + \norma{\nabla\phi}) \,.
  \non
\end{align}
\pier{Now, by combining it with \last{\eqref{pier12} and \eqref{pier13}}, then
applying the Young inequality and integrating over $(0,t)$, we arrive at  
\begin{align}
  & \frac12 \intQt |\nabla\mu|^2
  + \frac 12 \iO |\nabla\phi (t)|^2
  \leq 
  c \intQt \bigl(  |\nabla n|^2 + |\Delta n|^2 + |\phi|^2 + |\nabla \phi|^2 \bigr)\,, 
  \non
\end{align}
so that, in view of \eqref{1stimad} and \eqref{da1stimad},}
we conclude~that
\Beq
  \norma{\nabla\mu}_{\L2H}  + \pier{\norma{\nabla \phi}_{\L\infty H}} \leq c \, \norma u_{\L2H} \,.
  \non
\Eeq
From this, \pier{\eqref{1stimad}} and \eqref{da1stimad}\an{,} we deduce~that
\Beq
  \norma\mu_{\L2V} + \pier{\norma\phi_{\L\infty V\L2W}} \leq c \, \norma u_{\L2H} \,.
  \label{2stimad}
\Eeq

\step
Consequence

A comparison \an{argument} in \an{equation} \eqref{dprima} then \an{readily} yields that
\Beq
  \norma{\dt\phi}_{\L2\Vp} \leq c \, \norma u_{\L2H} \,.
  \label{da2stimad}
\Eeq

\step
Third estimate

\an{Next, we} take \an{an arbitrary} $v\in\juerg{\L2V}$, test \eqref{dterza} by $v$, and integrate over $(0,T)$.
We have~that
\begin{align}
  & \pier{{}\int_Q \dt a \, v {}}
  = - \intQ \nabla a \cdot \nabla v
  \non
  \\
  & \quad {} 
  + \chia \intQ \bigl( a \nabla\sig_1 + a_2 \nabla\sig \bigr) \cdot \nabla v
  + \intQ \bigl( a - (a_1+a_2)a + u \bigr) v \,.
  \non
\end{align}
Just some of the terms on the \rhs\ need some treatment.
The first one is the following:
\Beq
  \chia \intQ a \nabla\sig_1 \cdot \nabla v
  \leq \norma a_{\LQ4} \, \norma{\nabla\sig_1}_{\LQ4} \, \norma{\nabla v}_{\LQ2}
  \leq c \, \norma u_{\L2H} \, \norma v_{\L2V}\,, 
  \non
\Eeq
where the last inequality is due to the regularity of $\sig_1$ and~\eqref{1stimadbis}.
For analogous reasons we have~that
\Beq
  \chia \intQ a_2 \nabla\sig \cdot \nabla v
  \leq \norma{a_2}_{\LQ4} \, \norma{\nabla\sig}_{\LQ4} \, \norma{\nabla v}_{\LQ2}
  \leq c \, \norma u_{\L2H} \, \norma v_{\L2V}\,, 
  \non
\Eeq
as well as
\begin{align}
  & - \intQ (a_1+a_2) \, a \, v 
  \leq \norma{a_1+a_2}_{\L2{\Lx4}} \, \norma a_{\L\infty\Ldue} \, \norma v_{\L2{\Lx4}}
  \non
  \\
  & \leq c \, \norma{a_1+a_2}_{\L2V} \, \norma a_{\L\infty H} \, \norma v_{\L2V}
  \leq c \, \norma u_{\L2H} \, \norma v_{\L2V} \,.
  \non
\end{align}
Since the other terms on the \rhs\ can be estimated in a \an{straightforward} way,
we conclude that
\Beq
   \pier{{}\int_Q \dt a \, v \last{=} {}}\ioT \< \dt a(t) , v(t) > \, dt
  \leq c \, \norma u_{\L2H} \, \norma v_{\L2V} \,\an{,}
  \non
\Eeq
\an{whence, due to the arbitrariness of $v$, this entails} that
\Beq
  \norma{\dt a}_{\L2\Vp} \leq c \, \norma u_{\L2H} \,.
  \label{3stimad}
\Eeq
By recalling Remark~\ref{Remcontdep}, we see that the proof 
of the uniqueness part of Theorem~\ref{Wellposedness} and of \eqref{contdep} is complete.
\qed


\section{An auxiliary result}
\label{AUXILIARY}
\setcounter{equation}{0}

This section is devoted to \an{stating and proving} an auxiliary result that will be used twice in the sequel.
\an{In the following sections, we will repeatedly analyze systems related to the state system \Ipbl. Since \jurgen{these} share a very similar mathematical structure, we have decided to introduce an abstract result that encompasses both cases to be analyzed later. Theorem \ref{Wellposednessab}, proved below, will be used to \pier{show} the well-posedness of the linearized system and the \Frechet\ differentiability of the solution operator.}
We fix some $\ustar$ satisfying \eqref{HPu} and the corresponding solution $\soluzstar$ given by Theorem~\ref{Wellposedness},
and we recall at once that $\h'(\phistar)$ and $F''(\phistar)$ are bounded, 
since $\phistar$ satisfies the separation property~\eqref{separation}.
Moreover, we~fix 
\Beq
  g_1 \,,\, g_4 \,,\, g_5 \in \LQ2 , \quad
  g_2 \in \L2V
  \aand 
  \gg_3 \in (\LQ2)^2 \,\an{,}
  \label{HPg}
\Eeq
\an{and we}  notice that $\gg_3$ is \pier{a} vector-valued \an{function}.
Nevertheless, we often prefer to write~$g_3$ 
(i.e.,~we do not use the boldface character)
for uniformity.
Then, we look for the solution $\soluz$ to the problem stated below.
We remark that the notation $\soluz$ adopted in this section
\an{is unrelated to the original problem~\Pbl.}
We look for a quintuple $\soluz$ with the regularity properties
\begin{align}
  & \phi \in \calX_1 := \H1\Vp \cap \L\infty V \cap \L2W\,,
  \label{regphiab}
  \\
  & \mu \in \calX_2 := \L2V \,,
  \label{regmuab}
  \\
  & a \in \calX_3 := \H1\Vp \cap \L\infty H \cap \L2V\,,
  \label{regaab}
  \\
  & n \in \calX_4 := \H1H \cap \L\infty V \cap \L2W\,,
  \label{regnab}
  \\
  & \sig \in \pier{\calX_4} 
  \,,
  \label{regsigab}
\end{align}
\Accorpa\Regsoluzab regphiab regsigab
that solves the variational equations
\begin{align}
  & \< \dt\phi , v >
  + \iO \nabla\mu \cdot \nabla v
  - \chiphi \iO \nabla n \cdot \nabla v
  \non
  \\
  & = - m \iO \phi v
  + \iO \h'(\phistar) \phi \, v
  + \iO g_1 \, v \,,
  \label{primaab}
  \\
  \separa
  & \iO \nabla\phi \cdot \nabla v
  = - \iO F''(\phistar) \phi\, v
  + \iO \mu v
  + \iO g_2 \, v \,,
  \label{secondaab}
  \\
  \separa
  & \< \dt a , v > 
  + \iO \nabla a \cdot \nabla v
  - \chia \iO (a\nabla\sigstar + \astar\nabla\sig) \cdot \nabla v
  \non
  \\
  & = - \iO \gg_3 \cdot \nabla v
  + \iO (a - 2\astar a) v
  + \iO g_4 \, v \,,
  \label{terzaab}
  \\
  \separa
  & \iO \dt n \, v
  + \iO \nabla n \cdot \nabla v
  - \chiphi \iO \phi v
  = \iO (\cphi\phi + \cn n + \csig\sig) v \,,
  \label{quartaab}
  \\
  \separa
  & \iO \dt\sig \, v
  + \iO \nabla\sig \cdot \nabla v
  = \iO (-\sig + \chia a - a\sigstar - \astar\sig) v
  + \iO g_5 \, v \,,
  \label{quintaab}
\end{align}
for every $v\in V$ and \aet,
and satisfies the initial condition
\Beq
  (\phi,a,n,\sig)(0) = (0,0,0,0) 
  \quad \aeO \,.
  \label{cauchyab}
\Eeq
\Accorpa\Pblab primaab cauchyab

\pier{Let us introduce the space} \pier{(cf.~\Regsoluzab)
\Beq
\calX:= \calX_1\times\calX_2\times\calX_3\times\calX_4\times\pier{{}\calX_4}
\label{defcalX}
\Eeq
for the solutions to \Pblab.}

\Bthm
\label{Wellposednessab}
\juerg{Let the} assumptions of Theorem~\ref{Wellposedness} \juerg{and \eqref{HPg} be fulfilled, 
and assume that $\ustar$ satisfies \eqref{HPu} and $\soluzstar$ is} the corresponding solution.
Then the problem \Pblab\ has a unique solution $\soluz$ satisfying \Regsoluzab, and the estimate
\Beq
  \norma\soluz_{\pier{\calX}}
  \leq K_3 \, \Bigl(
    \somma i15 \norma{g_i}_{\LQ2}
    + \norma{\nabla g_2}_{\LQ2}
  \Bigr)
  \label{stabilityab}
\Eeq
holds true with a constant $K_3>0$ that depends only on 
$\Omega$, $T$, the structure of the original system, the initial data, and~$\umax$.
\Ethm

\Bdim
\an{To establish} existence, we just prove formal estimates for the solution\an{, but those computations do suggest} that the same estimates
\juerg{can be} performed on the solution to the $\an{k}$-dimensional system obtained from the Faedo--Galerkin scheme
constructed by using the first $\an{k}$ eigenfunctions of the Laplace operator with homogeneous boundary conditions.
\juerg{These bounds can then be used} to pass to the limit as $\an{k}$ tends to infinity and to construct a solution to the problem
satisfying \Regsoluzab\ and~\eqref{stabilityab}.

\step
First a priori estimate

We test the \an{above} equations \juerg{\eqref{primaab}--\eqref{quintaab}}\an{, in the order,} by $\phi$, $-\Delta\phi$, $a$, $\dt n$ and $\dt\sig-\Delta\sig$, 
respectively.
We obtain~that
\begin{align}
  & \frac 12 \, \ddt \iO |\juerg{\phi}|^2
  + \iO \nabla\mu \cdot \nabla\phi
  + m \iO |\phi|^2
  \non
  \\
  & = \chiphi \iO \nabla n \cdot \nabla\phi
  + \iO \h'(\phistar) |\phi|^2
  + \iO g_1 \, \phi \,,
  \non
  \\
  \separa
  & \iO |\Delta\phi|^2
  = \iO F''(\phistar) \phi\, \Delta\phi
  + \iO \nabla\mu \cdot \nabla\an{\phi}
  - \iO g_2 \, \Delta\phi \,,
  \non
  \\
  \separa
  & \frac 12 \, \ddt \iO |a|^2
  + \iO |\nabla a|^2
  - \chia \iO (a\nabla\sigstar + \astar\sig) \cdot \nabla a
  \non
  \\
  & = \iO |a|^2
  - 2 \iO \astar |a|^2
  - \iO \gg_3 \cdot \nabla a
  + \iO g_4 \, a\,,
  \non
  \\
  \separa
  & \iO |\dt n|^2
  + \frac 12 \, \ddt \iO |\nabla n|^2
  = \chiphi \iO \phi \, \dt n
  + \iO (\cphi\phi + \cn n + \csig\sig) \, \dt n \,,
  \non
  \\
  \separa
  & \iO |\dt\sig|^2
  + \an{\, \ddt \iO |\nabla\sig|^2}
  + \iO |\Delta\sig|^2
  \non
  \\
  & = \iO (-\sig + \chia a - a\sigstar - \astar\sig) (\dt\sig - \Delta\sig)
  + \iO g_5 \, (\dt\sig - \Delta\sig)\,.
  \non
\end{align}
\an{Then,} we take the sum of these identities \an{and} add $(1/2)d/dt\iO|n|^2$ and $d/dt\iO|\sig|^2$ to the \lhs\ of the resulting equality 
and the same terms, written in the form $\iO n\dt n$ and $2\iO\sig\dt\sig$, to \an{the} \rhs.
\juerg{The \lhs\ then becomes}
\begin{align} 
  & \frac 12 \, \ddt \iO |\phi|^2
  + m \iO |\phi|^2
  + \iO |\Delta\phi|^2
  + \frac 12 \ddt \iO |a|^2 
  + \iO |\nabla a|^2
  \non
  \\
  & \quad {}
  + \iO |\dt n|^2
  + \frac 12 \, \ddt \iO \bigl(\an{|n|^2+ |\nabla n|^2 }\bigr)
  \non
  \\
  & \quad {}
  + \iO |\dt\sig|^2
  + \ddt \iO \bigl(\an{|\sig|^2 + |\nabla\sig|^2 }\bigr)
  + \iO |\Delta\sig|^2\,,
  \non
\end{align}
and we have to estimate the terms on the corresponding \rhs.
However, two of these cancel each other, one is nonpositive,
and most of the others can be easily dealt with using Young's inequality.
So, we \juerg{discuss} only the most delicate \an{terms}.
As usual, we \juerg{intend to apply Gronwall's} lemma after time integration.
The first term we consider is the following:
\begin{align}
  & \chia \iO a \nabla\sigstar \cdot \nabla a
  \leq \norma a_2 \, \norma{\nabla\sigstar}_\infty \, \norma{\nabla a}_2
  \leq \frac 14 \iO |\nabla a|^2
  + \norma{\nabla\sigstar}_\infty^2 \iO |a|^2 \,.
  \non
\end{align}
To deal with the next one, we owe to \juerg{the second inequality in \eqref{lady} to find} that
\begin{align}
  & \chia \iO \astar \nabla\sig \cdot \nabla a
  \leq \norma\astar_4 \, \norma{\nabla\sig}_4 \, \norma{\nabla a}_2
  \leq \frac 14 \iO |\nabla a|^2
  + \norma\astar_4^2 \, \norma{\nabla\sig}_4^2
  \non
  \\
  & \leq \frac 14 \iO |\nabla a|^2
  + \norma\astar_4^2 \, \normaV\sig \, (\norma\sig + \norma{\Delta\sig})
  \non
  \\
  & \leq \frac 14 \iO |\nabla a|^2
  + \frac 14 \iO |\Delta\sig|^2
  + c \, \an{\bigl(\norma\astar_4^2+\norma\astar_4^4 \bigr)} \, \normaV\sig^2 
  \,.
  \non
\end{align}
We notice that the functions $\,\,t\mapsto\norma{\nabla\sigstar(t)}_\infty^2$ \an{and} $t\mapsto \an{(\norma{\astar(t)}_4^2 + \norma{\astar(t)}_4^4)}$
belong to $L^1(0,T)$.
The last integral we consider is the following:
\begin{align}
  & \iO \astar\sig (\dt\sig - \Delta\sig)
  \leq \norma\astar_4 \, \norma\sig_4 \, (\norma{\dt\sig}_2 + \norma{\Delta\sig}_2)
  \non
  \\
  & \leq \frac 18 \iO \bigl( |\dt\sig|^2 + |\Delta\sig|^2 \bigr)
  + c \, \norma\astar_4^2 \, \normaV\sig^2 \,.
  \non
\end{align}
By collecting, rearranging, and applying the Gronwall lemma, we conclude that
\begin{align}
  & \norma\phi_{\L\infty H\cap\L2W}
  + \norma a_{\L\infty H\cap\L2V}
  + \norma n_{\H1H\cap\L\infty V}
  \non
  \\
  & \quad {}
  + \norma\sig_{\H1H\cap\L\infty V\cap\L2W}
  \leq c \, \somma i15 \norma{g_i}_\juerg{L^2(Q)} \,.
  \label{1stimaab}
\end{align}

\step
Consequence

Notice that \eqref{1stimaab} also yields an estimate for $a$ and $\nabla\sig$ in $\LQ4$, thanks to  \pier{\eqref{disug}.}
Moreover, by comparison, first in \eqref{quartaab} and then in \eqref{secondaab}, \pier{and elliptic regularity}, we derive estimates for $n$ and~$\mu$.
In conclusion, we have~that
\begin{align}
  & \norma a_{\LQ4}
  + \norma{\nabla\sig}_{\LQ4}
  + \norma n_{\L2W}
  + \norma\mu_{\L2H}
  \leq c \, \somma i15 \norma{g_i}_\juerg{L^2(Q)} \,.
  \label{da1stimaab}
\end{align}

\step
Second a priori estimate

We test \eqref{terzaab} by an arbitrary $v\in\L2V$ and integrate over~$(0,T)$ \an{to} obtain~that
\begin{align}
  & \ioT \< \dt a(t) , v(t) > \, dt
  = - \intQ \nabla a \cdot \nabla v
  + \chia \intQ (a \nabla\sigstar + \astar\nabla\sig) \cdot \nabla v
  \non
  \\
  & \quad {}
  - \intQ \gg_3 \cdot \nabla v
  + \intQ (a + g_4) v
  - \intQ 2\astar a v \,.
  \non
\end{align}
Just some of the terms on the \rhs\ need a treatment.
We have~that
\begin{align}
  & \chia \intQ (a \nabla\sigstar + \astar\nabla\sig) \cdot \nabla v
  \non
  \\
  & \leq \norma a_{\LQ4} \, \norma{\nabla\sigstar}_{\LQ4} \, \norma{\nabla v}_{\LQ2}
  + \norma\astar_{\LQ4} \, \norma{\nabla\sig}_{\LQ4} \, \norma{\nabla v}_{\LQ2}
  \non
  \\
  & \leq c \, \norma v_{\L2V} \, \somma i15 \norma{g_i}_\juerg{L^2(Q)}\,, 
  \non
\end{align}
and we similarly obtain that
\Beq
  - \intQ 2\astar a v
  \leq c \, \norma v_{\L2V} \, \somma i15 \norma{g_i}_\juerg{L^2(Q)} \,.
  \non
\Eeq
Hence, \pier{it turns out that} 
\Beq
  \ioT \< \dt a(t) , v(t) > \, dt
  \leq c \, \norma v_{\L2V} \, \somma i15 \norma{g_i}_\juerg{L^2(Q)} \,\an{,}
  \non
\Eeq
\an{and, s}ince $v$ is arbitrary in $\L2V$, this means that
\Beq
  \norma{\dt a}_{\L2\Vp}
  \leq c \somma i15 \norma{g_i}_\juerg{L^2(Q)} \,.
  \label{2stimaab}
\Eeq

\step
Third a priori estimate 

We test \eqref{primaab} by both $\mu$ and $g_2-F''(\phistar)\phi$.
At the same time, we test \eqref{secondaab} by~$\dt\phi$.
We obtain~that
\begin{align}
  & \iO \dt\phi \, \mu 
  + \iO |\nabla\mu|^2
  = \chiphi \iO \nabla n \cdot \nabla\mu
  - m \iO \phi \mu
  + \iO \h'(\phistar) \, \phi \mu
  + \iO g_1 \mu \,,
  \non
  \\
  & \iO \dt\phi \, \bigl( g_2-F''(\phistar)\phi \bigr)
  = \iO \bigl( -\nabla\mu + \chiphi \nabla n \bigr) \cdot \nabla \bigl( g_2-F''(\phistar)\phi \bigr)
  \non
  \\
  & \quad {}
  + \iO \bigl( -m\phi + \h'(\phistar)\juerg{\phi} + g_1 \bigr) \bigl( g_2-F''(\phistar)\phi \bigr)\,, 
  \non
  \\
  & \frac 12 \, \ddt \iO |\nabla\phi|^2
  = \iO \mu \, \dt\phi
  + \iO \bigl( g_2 - F''(\phistar) \, \phi \bigr) \dt\phi\,.
  \non
\end{align}
At this point, we add these equalities to each other and notice several cancellations.
The \lhs\ \juerg{then becomes}
\Beq
  \iO |\nabla\mu|^2
  + \frac 12 \, \ddt \iO |\nabla\phi|^2\,,
  \label{lhsterzaab}
\Eeq
and we \juerg{need to} estimate the terms on the corresponding \rhs.
However, we only treat the most delicate of them,
since the others are simple \juerg{to handle}.
By \juerg{recalling}, in particular, that $\nabla\phistar$ is bounded \an{owing to the regularity in \eqref{nablaphibdd}}, we have~that
\begin{align}
  &
  \iO \bigl( \juerg{-}\nabla\mu + \chiphi \nabla n \bigr) \cdot \nabla \bigl( \last{g_2 -} F''(\phistar)\phi \bigr)
  \non
  \\
  & \leq \frac 14 \iO |\nabla\mu|^2
  + c \iO |\nabla n|^2
  \last{+ c \iO |\nabla g_2|^2} 
  + c \iO |\phi \, F'''(\phistar) \nabla\phistar|^2 
  + c \iO |F''(\phistar) \nabla\phi|^2
  \non
  \\
  & \leq \frac 14 \iO |\nabla\mu|^2
  + c \iO |\nabla n|^2
   \last{+ c \iO |\nabla g_2|^2} 
 + c \iO |\phi|^2 
  + c \iO |\nabla\phi|^2\,,
  \non
\end{align}
and \an{we} account for \eqref{1stimaab}.
By combining it with the \juerg{already proved estimates},
we conclude~that
\Beq
  \norma\mu_{\L2V} + \norma\phi_{\L\infty V}
  \leq c \somma i15 \norma{g_i}_{\LQ2}
  + c \, \norma{\nabla g_2}_{\LQ2}\,.
  \label{3stimaab}
\Eeq

\step
Conclusion
 
\juerg{Now that} $\mu$ is estimated in $\L2V$,
a comparison \an{argument} in \eqref{primaab}\an{, using also \eqref{3stimaab},} yields a similar estimate for $\dt\phi$ in $\L2\Vp$,
and the existence part of the proof is complete.

\step
Uniqueness

By linearity, we consider only the homogeneous problem.
We want to come back to the proof of \eqref{1stimaab} 
and show that the procedure used there can be made rigorous.
On account of the regularity \Regsoluzab,
we notice that the equations \eqref{secondaab}, \eqref{quartaab}, and \eqref{quintaab},
can be written in \last{strong} form.
For instance, \eqref{secondaab} with $g_2=0$ becomes
\Beq
  -\Delta\phi = -F''(\phistar) \phi + \mu 
  \quad \aeQ \,.
  \label{strongsecondaab}
\Eeq
Hence, instead of testing \eqref{secondaab} by $-\Delta\phi$, 
we can multiply \eqref{strongsecondaab} by $-\Delta\phi$ and integrate over~$\Omega$.
Therefore, we can still arrive at \eqref{1stimaab}, which yields $(\phi,a,n,\sig)=(0,0,0,0)$\last{, from which}, \eqref{strongsecondaab} implies that $\mu=0$ as well.
\Edim


\section{The control problem}
\label{CONTROL}
\setcounter{equation}{0}

In this section, we give the first result on the control problem \an{presented} in the Introduction.
For the reader's convenience, we recall the definitions of the cost functional $\calJ$
and of the set $\Uad$ of the admissible controls:
\begin{align}
  & \J (\phi, u)
  := \frac {\b1}2 \intQ |\phi - \phiQ|^2
  + \frac {\b2}2  \iO |\phi(T) - \phiO|^2
  + \frac {\b3}2 \intQ |u|^2 
  \non
  \\
  & \quad \hbox{for $\phi\in\last{\C0H}$ and $u\in\LQ2$},
  \label{cost}
  \\[2mm]
  & \Uad : = \big\{ u \in \calU : 0 \leq u \leq u_{\rm max} \ \aeQ \big\}, 
  \quad \hbox{where} \quad
  \calU := \LQ\infty .
  \label{Uad}
\end{align}
\juerg{We make the following assumptions:}
\begin{align}
  & \b i \in [0,+\infty) \, \hbox{ for $i=1,2,3$}\pier{, \hbox{ with } \, b_3>0;}  \quad
  \phiQ \in \LQ2, \, \   \phiO \in V .
  \label{HPcost}
  \\
  & \hbox{$\umax\in\LQ\infty$ \ is nonnegative} .
  \label{HPumaxbis}
\end{align}
Then the control problem is given by:
\begin{align}
  & \hbox{Minimize} \ \ \J (\phi, u) \ \ \text{subject to \,\,$u\in\Uad$\,\, and to the constraint that}
  \non
  \\
  & \mbox{$\soluz$ is the solution to the system \Pbl}.
  \label{control}
\end{align}   
In the remainder of the paper, it is understood that the above assumptions are in force,
as well as those on the structure and the data (with the same $\umax$ as here, of course)
that ensure well-posedness for the state system (see Theorem~\ref{Wellposedness}).
\juerg{We therefore do not recall them in any of the following statements.}
\an{Besides, by virtue of Theorem \ref{Wellposedness}, we can introduce the {\it control-to-state} operator $\calS$ as}
\begin{align*}
\pier{\calS\juerg{{}:=(\calS_1,\calS_2,\calS_3,\calS_4,\calS_5)} \, \hbox{ mapping } \, u \in \Uad \, \hbox{ into }\,  \soluz \in \pier{\calY}},
\end{align*}
where \pier{$\calY$ is the regularity space defined in \eqref{regtot}. The same operator is \Lip\ continuous from $\Uad$ into the space $\calX$ specified by \eqref{defcalX}, in the sense of the continuous dependence estimate~\eqref{contdep}.}

\Bthm
\label{Optimum}
The control problem \eqref{control} has at least one solution, \last{that is},
there exists at least one $\ustar\in\Uad$ that satisfies
\Beq
  \calJ(\phistar,\ustar) \leq \calJ(\phi,u)
  \quad \hbox{for every $u\in\Uad$}\pier{,}
  \label{optimum}
\Eeq
where $\phistar$ and $\phi$ are the first components of the solutions to the state system \Pbl\
corresponding to $\ustar$ and~$u$, respectively.
\Ethm

\Bdim
We use the direct method \an{of the calculus of variations}.
\an{To begin with, we recall \eqref{cost} and observe that \pier{$\J (\phi, u)$} 
is bounded from below as it is \jurgen{nonnegative}.}
Thus, we term $\Lambda$ the infimum of the cost functional under the constraints of the control problem
and fix a minimizing sequence $\graffe{\uk}$ and the sequence of the corresponding solutions \an{$\soluzk\in \calY$ as given by Theorem \ref{Wellposedness}}.
\an{Namely}, we have~that
\Beq
  \lim_{k\to\infty} \calJ(\phik,\uk) = \Lambda \,.
  \non
\Eeq
Now, $\Uad$ is bounded in $\LQ\infty$, and all of the solutions $\soluzk$ satisfy the stability estimate \eqref{stability}.
Hence, we can use well-known compactness results \juerg{to} obtain~that\an{, possibly for a nonrelabeled subsequence, as $k\to\infty$,}
\an{
\begin{align}
 &  \uk \to \ustar
  \quad \hbox{weakly star in $\LQ\infty$},
  \label{convu}
  \\
  & (\phik , \muk, \ak, \nk, \sigk) \to \pier{(\phistar , \mustar, \astar , \nstar, \sigstar)}
  \quad \hbox{weakly star in $\calY$},
  \label{convall}
\end{align}}%
\pier{for some limiting functions $\ustar$ and $\phistar , \mustar, \astar , \nstar, \sigstar$.}
\juerg{At this point}, strong convergence properties are needed,
and we apply \cite[Sect.~8, Cor.~4]{Simon} several times.
First, we notice that $\graffe{\phik}$ \pier{weakly star converges in} \an{$\calY_1$ that is compactly embedded in~$\an{\CQ0}$.}
We \an{thus} infer~that\an{, as $k\to\infty$,}
\begin{align}
  & \phik \to \phi
  \quad \hbox{strongly in \pier{$\CQ0$}},
  \quad \hbox{whence also}
  \non
  \\
  & F'(\phik) \to F'(\phi)
  \aand
  \h(\phik) \to \h(\phistar)
  \quad \hbox{strongly in $\CQ0$},
  \non
\end{align}
since the functions $\phik$ satisfy the separation property \eqref{separation}
and $F'$ and $\h$ are \Lip\ continuous in the interval $[r_-,r_+]$.
Moreover, \an{\eqref{convall} (cf. \eqref{rega})} implies that $\ak$ converges to $\astar$ strongly in $\an{\C0{H} \cap \L2{V}}$. \pier{In view of \eqref{disug},
we deduce~that\an{, as $k\to\infty$,} $\ak \to \astar$ strongly in $\LQ4$ and consequently 
\Beq
  (\ak)^2 \to (\astar)^2
  \quad \hbox{strongly in $\LQ2$} \,.
  \non
\Eeq
Next, from \eqref{convall} and \eqref{regsig} it follows that\an{, as $k\to\infty$,}}
\Beq
  \sigk \to \sigstar
  \quad \pier{\hbox{strongly in $C^0(\ov Q) \cap \L4{\Wx{1,4}}$},}
  \non
\Eeq
\pier{whence
\Beq
  \ak\sigk \to \astar\sigstar
  \aand
  \ak \nabla\sigk \to \astar\nabla\sigstar
  \quad \hbox{strongly in ${\LQ2}$ and $\LQ2^2$}.
  \non
\Eeq
Collecting all this information, and passing to the limit in the variational equalities \eqref{prima}--\eqref{quinta} written for $(\phik , \muk, \ak, \nk, \sigk)$}, we deduce that
$\soluzstar$ solves the state system corresponding to~$\ustar$.
This shows that $\soluzstar$ is \pier{actually $\calS(\ustar)$.}
\an{Thus, we have} that
\Beq
  \calJ(\phistar,\ustar)
  \leq \liminf_{k\to\infty} \calJ(\phik,\uk)
  = \Lambda\,,
  \non
\Eeq 
so that $\calJ(\phistar,\ustar)=\Lambda$ and $\ustar$ is an optimal control.
\Edim


\section{Necessary conditions for optimality}
\label{NECESSARY}
\setcounter{equation}{0}

\an{T}o obtain significant necessary optimality conditions,
\an{we} prove some  differentiability \an{property of  the control-to-state operator $\calS$.} 
\juerg{To this end, recall} the definitions of the spaces $\calY_i$ and $\calX_i$, $i\pier{{}=1,...,4{}}$,
given in \Regsoluz\ and \Regsoluzab, respectively\juerg{. Notice that} $\calY_i\emb\calX_i$ for every~$i\pier{{}=1,...,4{}}$.
\pier{We then recall the definition~\eqref{defcalX} of $\calX$ and consider the mapping} $\calS:\Uad\to\calX$ by observing that 
\begin{align}
  & \hbox{for $u\in\Uad$, \ \juerg{$\calS(u) =\soluz $ }}\, 
	\mbox{is the solution to \Pbl.}
	\label{defS}
\end{align}

\an{Often, it is possible to extend $\calS$ to an open subset of $\calU$ and prove the differentiability of this extension.}
However, we cannot develop this idea \an{in our scenario}.
Indeed, although the constraint $u\leq\umax$ could be replaced by $\norma u_\infty<K$
(with any prescribed $K>\norma\umax_\infty$)  without any significant change in our previous proofs,
we cannot avoid the constraint $u\geq0$.
Therefore, we can only prove some kind of \an{\emph{sectorial} differentiability}
that is close to \Frechet\ differentiability but \an{intrinsically} involves a constraint for the increments.
Also in the present case, a crucial role is played by the linearized system we introduce at once.
To this end, we fix $\ustar\in\Uad$ and $\soluzstar:=\calS(\ustar)$\an{.}
By accounting for the separation property \eqref{separation} 
and for the smoothness of $F$ on the interval $[r_-,r_+]$,
we infer~that
\Beq
  F^{(j)}(\phistar) \in \juerg{L^\infty(Q)}
  \aand
  \norma{F^{(j)}(\phistar)}_\infty \leq K_1'
  \quad \hbox{for $\last{j\in \{0,...,4\}}$},
  \label{bddFj}
\Eeq
where $K_1'$ is similar to the constant $K_1$ appearing in~\eqref{stability}.

Let us come to the linearized problem associated with~$\ustar$.
For \an{a given} $h\in\LQ2$, it consists in looking for a quintuple $\soluzl$ 
with the regularity
\Beq
  \soluzl \in \calX
  \label{regsoluzl}
\Eeq
that solves the variational equations
\begin{align}
  & \< \dt\psi , v >
  + \iO \nabla\eta \cdot \nabla v
  - \chiphi \iO \nabla\nu \cdot \nabla v
  \non
  \\
  & = - m \iO \psi \, v
  + \iO \h'(\phistar) \, \psi v \,,
  \label{primal}
  \\
  & \iO \nabla\psi \cdot \nabla v
  + \iO F''(\phistar) \, \psi v 
  = \iO \eta v\,,
  \label{secondal}
  \\
  & \< \dt\al , v >
  + \iO \nabla\al \cdot \nabla v
  - \chia \iO (\al \nabla\sigstar + \astar \nabla\omega) \cdot \nabla v
  \non
  \\
  & = \iO (\al - 2\astar\al) v
  + \iO h v\,,
  \label{terzal}
  \\
  & \iO \dt\nu \, v
  + \iO \nabla\nu \cdot \nabla v
  - \chiphi \iO \psi v
  = \iO (\cphi\psi + \cn\nu + \csig\omega) v\,,
  \label{quartal}
  \\
  & \iO \dt\omega \, v
  + \iO \nabla\omega \cdot \nabla v
  = \iO (-\omega + \chia\al - \al\sigstar - \astar\omega) v\,,
  \label{quintal}
\end{align}
\juerg{for} every $v\in V$ and \aet,
and satisfies the initial condition
\Beq
  (\psi,\al,\nu,\omega)(0) = (0,0,0,0) \,.
  \label{cauchyl}
\Eeq
\Accorpa\Pbll primal cauchyl

We have the following result \an{concerning well-posedness.}

\Bthm
\label{Wellposednessl}
Let $\ustar\in\Uad$ be given and $\soluzstar:=\calS(\ustar)$.
Then, for every $h\in\LQ2$, problem \Pbll\ has a unique solution $\soluzl$ 
satisfying \eqref{regsoluzl}, and the estimate
\Beq
  \norma\soluzl_\calX \leq K_3 \, \norma h_{\LQ2}
  \label{stimal}
\Eeq 
holds true with a constant $K_3>0$ 
that depends only on $\Omega$, $T$, the structure \pier{and the data of the state system}, and~$\umax$.
\Ethm

\Bdim
It suffices to apply Theorem~\ref{Wellposednessab}.
Indeed, problem \Pbll\ is the particular case of problem \Pblab, 
where $g_4=h$ and \an{$g_1=g_2=g_3=g_5=0$}.
\Edim

Let us \an{move} to \an{prove} some differentiability property of \an{the solution operator}~$\calS$.
As said before, we cannot speak of \Frechet\ differentiability.
However, the above result implies that the linear mapping $\,h\mapsto\soluzl\,$ is continuous from $\LQ2$ into $\calX$,
and we are now going to see that it plays \juerg{a similar role as a} \Frechet\ derivative.
Indeed, we can prove the following result\an{.}

\Bthm
\label{Frechet}
Let $\ustar\in\Uad$ and $\soluzstar:=\calS(\ustar)$. 
For every $h\in\LQ2$, let $\soluzl$ be 
\juerg{the} solution to the corresponding linearized system \Pbll.
Then we have that
\begin{align}
  & \frac {\norma{\calS(\ustar+h) - \calS(\ustar) - \soluzl}_{\calX}} {\norma h_{\LQ2}}
  \quad \hbox{tends to zero}
  \non
  \\
  & \quad \hbox{as \ $\norma h_{\LQ2}$ \ tends to zero under the constraint that \ $\ustar+h\in\Uad$}.
  \label{frechet}
\end{align}
\Ethm

\Bdim
We assume that $\ustar+h$ belongs to~$\Uad$
and introduce $\soluzh:=\calS(\ustar+h)$ and the quintuplet $\soluzF\in\calX$ defined~by
\begin{align}
  & \ph := \phih-\phistar-\psi \,, \quad
  \rho := \muh-\mustar-\eta \,, \quad
  \gamma := \ah-\astar-\al\,, 
  \non
  \\
  & \quad \lam := \nh-\nstar-\nu 
  \aand
  \xi := \sigh-\sigstar-\omega \,,
  \label{soluzF}
\end{align}
so that
\Beq
  \calS(\ustar+h) - \calS(\ustar) - \soluzl
  = \soluzF \,.
  \non
\Eeq
Then, $\soluzF$ satisfies the variational equations
\begin{align}
  & \< \dt\ph , v >
  + \iO \nabla\rho \cdot \nabla v
  - \chiphi \iO \nabla\lam \cdot \nabla v
  \non
  \\
  & = - m \iO \ph v
  + \iO \an{\big[}\h(\phih) - \h(\phistar) - \h'(\phistar) \psi \an{\big]} \, v\,,
  \label{primaF}
  \\
  & \iO \nabla\ph \cdot \nabla v
  = \iO \rho v
  - \iO \an{\big[} F'(\phih) - F'(\phistar) - F''(\phistar) \psi \an{\big]} \, v\,,
  \label{secondaF}
  \\
  & \< \dt\juerg{\gamma} , v > 
  + \iO \nabla\juerg{\gamma} \cdot \nabla v
  - \chia \iO \an{\big[}\ah\nabla\sigh - \astar\nabla\sigstar - \al\nabla\sigstar - \astar\nabla\omega\an{\big]} \cdot \nabla v
  \non
  \\
  & = \iO \gamma v
  - \iO \bigl[ (\ah)^2 - (\astar)^2 - 2\astar\al \bigr] \, v\,,
  \label{terzaF}
  \\
  & \iO \dt\lam \, v
  + \iO \nabla\lam \cdot \nabla v
  - \chiphi \iO \ph v
  = \iO (\cphi\ph + \cn\gamma + \csig\xi) v\,,
  \label{quartaF}
  \\
  & \iO \dt\xi \, v
  + \iO \nabla\xi \cdot \nabla v
  = \iO \bigl( -\xi + \chia\gamma - \bigl[\ah\sigh-\astar\sigstar-\al\,\juerg{\sigstar}-\astar\omega\bigr] \bigr) v\,,
  \label{quintaF}
\end{align}
all for every $v\in V$ and \aet,
and the initial condition
\Beq
  (\ph,\gamma,\lam,\an{\xi})(0) = (0,0,0,0) 
  \quad \aeO \,.
  \label{cauchyF}
\Eeq
\Accorpa\PblF primaF cauchyF
We transform some terms on the \rhs s, \juerg{using Taylor expansions in three cases}.
We have~that
\begin{align}
  & \h(\phih) - \h(\phistar) - \h'(\phistar) \psi
  = \h'(\phistar) \, \ph
  + R_1 (\phih-\phistar)^2
  \non
  \\
  & \quad \hbox{with} \quad
  R_1 = \int_0^1 (1-s) \h''(\phistar+s(\phih-\phistar)) \, ds\,,
  \non
  \\
  & F'(\phih) - F'(\phistar) - F''(\phistar) \psi
  = F''(\phistar) \, \ph
  + R_2 (\phih-\phistar)^2
  \non
  \\
  & \quad \hbox{with} \quad
  R_2 = \int_0^1 (1-s) F'''(\phistar+s(\phih-\phistar)) \, ds\,,
  \non
  \\
  & \ah\nabla\sigh - \astar\nabla\sigstar - \al\nabla\sigstar - \astar\nabla\omega
  = (\ah-\astar) (\nabla\sigh-\nabla\sigstar) + \gamma \nabla\sigstar + \astar \nabla\xi\,,
  \non
  \\[2mm]
  & (\ah)^2 - (\astar)^2 - 2\astar\al
  = 2\astar\gamma + R_3 (\ah-\astar)^2 
  \non
  \\
  & \quad \hbox{with} \quad
  R_3 = \int_0^1 (1-s) \, 2(\phistar+s(\phih-\phistar)) \, ds\,,
  \non
	\\
	&\juerg{\ah\sigh-\astar\sigstar-\al\sigstar-\astar\omega
	=(\ah-\astar)\,(\sigh-\sigstar) + \gamma\sigstar +\astar \xi\,.}
\end{align}
We notice that both $\phistar$ and $\phih$ satisfy the separation condition \eqref{separation},
so that the same holds for $\phistar+s(\phih-\phistar)$ for every $s\in[0,1]$.
Hence, the quantities under the above integrals over $(0,1)$ are bounded,
and all the \an{remainders $R_1, ...,R_\last{3}$} are \an{uniformly bounded}.
At this point, we notice that problem \PblF\ takes the form \Pblab\ with
\begin{align}
  & g_1 = R_1 \, (\phih-\phistar)^2 \,, \quad
  g_2 = \pier{-{}} R_2 \, (\phih-\phistar)^2 \,, \quad
  \gg_3 = -\chia (\ah-\astar) \, \nabla(\sigh-\sigstar)\,,
  \non
  \\
  & g_4 =\pier{-{}} R_3 \, (\ah-\astar)^2
  \aand
  g_5 = \pier{-{}} \juerg{(\ah-\astar)} \, (\sigh-\sigstar) \,.
  \non
\end{align}
Hence, we can apply Theorem~\ref{Wellposednessab}, \an{to infer} that
\Beq
  \norma\soluzF_\calX
  \leq K_3 \, \Bigl(
    \somma i15 \norma{g_i}_{\LQ2}
    + \norma{\nabla g_2}_{\LQ2}
  \Bigr) \,.
  \label{stabilityF}
\Eeq
Thus, it remains to estimate the \rhs\ of \pier{\eqref{stabilityF}}.
In doing this, we also account for Theorem~\ref{Contdep}
and apply \eqref{contdep} to $\soluzh$ and $\soluzstar$,
possibly combined with the last inequality in \eqref{disug}\juerg{. We} have~that
\begin{align}
  & \norma{g_1}^2_{\LQ2}
  + \norma{g_2}^2_{\LQ2}
  \leq c \intQ |\phih-\phistar|^4
  \leq c \, \norma h_{\LQ2}^4\,,
  \non
  \\
  & \norma{\gg_3}^2_{\an{\LQ2}}
  \leq \iO |\ah-\astar|^2 \, |\nabla(\sigh-\sigstar)|^2
  \non
  \\
  &\quad{}
  \leq \norma{\ah-\astar}_{\LQ4}^2 \, \norma{\nabla(\sigh-\sigstar)}_{\LQ4}^2
  \leq c \, \norma h_{\LQ2}^4\,,
  \non
  \\
  & \norma{g_4}^2_{\LQ2}
  \leq c \intQ |\ah-\astar|^4
  \leq c \, \norma h_{\LQ2}^4\,,
  \non
  \\
  & \norma{g_5}^2_{\an{\LQ2}}
  \leq \iO |\ah-\astar|^2 \, |\sigh-\sigstar|^2
  \non
  \\
  &\quad{}
  \leq \norma{\ah-\astar}_{\LQ4}^2 \, \norma{\sigh-\sigstar}_{\LQ4}^2
  \leq c \, \norma h_{\LQ2}^4 \,.
  \non
\end{align}
Finally, since it turns out that $|\nabla R_2|\leq c(|\nabla\phistar|+|\nabla\phih|)\leq c$,
we also have~that
\begin{align}
  & \norma{\nabla g_2}_{\LQ2}^2
  = \norma{(\phih-\phistar)^2\,\nabla R_2+2R_2(\phih-\phistar)\nabla(\phih-\phistar)}_{\LQ2}^2
  \non
  \\
  & \leq c \, \bigl( \norma{\phih-\phistar}_{\LQ4}^4 + \norma{\phih-\phistar}_{\LQ4}^2 \, \norma{\nabla(\phih-\phistar)}_{\LQ4}^2
  \leq c \, \norma h_{\LQ2}^4 \,.
  \non
\end{align}
Therefore, \eqref{stabilityF} implies that
\Beq
  \norma\soluzF_\calX
  \leq c \, \norma h_{\LQ2}^2\,,
  \non
\Eeq
and \eqref{frechet} \an{readily} follows.
\Edim

\juerg{Thanks to the above} result and the convexity of~$\Uad$, \juerg{a standard argument leads to the following
necessary condition for an element $\ustar\in\Uad$ to be an optimal control:}
\begin{align}
  & \b1 \intQ (\phistar-\phiQ) \, \psi
  + \b2 \iO \bigl( \phistar(T)-\phiO \bigr) \, \psi(T)
  + \b3 \intQ \ustar \, (u-\ustar)
  \geq 0
  \non
  \\
  & \quad \hbox{for every $u\in\Uad$}\,,
  \label{badnc}
\end{align}
where $\psi$ is the first component of the solution to the linearized problem \Pbll\ corresponding to $h:=u-\ustar$.
However, this condition is \an{problematic}, since it requires to solve the linearized problem infinitely many times
because $u$ is arbitrary in~$\Uad$.
As usual, this trouble is overcome by introducing a proper adjoint problem \juerg{associated with a given $\ustar\in\Uad$.
In order to simplify its presentation, we use some abbreviations:
for some pairs $(i,j)$, \an{with $i,j \in \{\last{0},...,5\}$,}} we define $\last{f_{i,j}}$ as follows: \an{
\begin{align}
  & \f10 = \b1 (\phistar-\phiQ)
 	\,,\ 
  \label{deff10}                                                           
  \f11 = m-\h'(\phistar) \,,\ 
  \f12 = F''(\phistar) \,,\ 
  \f14 = -\chiphi-\cphi \,,\ 
  \\
  & 
  \f33 = -1 \,,\ 
  \f35 = \sigstar-\chia\,,
  \label{deff3j}
  \\
  & 
  \f54 = -\csig \,,\ 
  \f55 = 1 \,\an{,}
  \label{deff5j}
\end{align}
and put to zero all the other cases.}
We notice that $\f10\in\LQ2$ and that every other $\f ij$ is a bounded function.
At this point, we can write the adjoint problem associated with~$\ustar$.
It consists in looking for a quintuple $\soluzad$ with the regularity properties
\begin{align}
  & \p1 \in \H1\Vp \cap \L\infty V \cap \L2{\an{W\cap \Hx3}}\,,
  \label{regprimaad}
  \\
  & \p2 \in \L2V\,,
  \label{regsecondaad}
  \\
  & \p3 \in \H1H \cap \L\infty V \cap \L2W\,,
  \label{regterzaad}
  \\
  & \p4 \in \H1V \cap \L2{\an{W\cap \Hx3}}\,,
  \label{regquartaad}
  \\
  & \p5 \in \H1\Vp \cap \L\infty H \cap \L2V\,,
  \label{regquintaad}
\end{align}
\Accorpa\Regsoluzad regprimaad regquintaad
that solves the variational equations
\begin{align}
  & - \< \dt\p1 , v >
  + \iO \nabla\p2 \cdot \nabla v
  + \somma j15 \iO \f1j \, \p j \, v
  = \iO \f10 \, v\,,
  \label{primaad}
  \\
  & \iO \nabla\p1 \cdot \nabla v
  = \iO \p2 \, v\,,
  \label{secondaad}
  \\[2mm]
  & - \iO \dt\p3 \, v 
  + \iO \nabla\p3 \cdot \nabla v
  - \chia \iO (\nabla\sigstar \cdot \nabla\p3) v
  \non
  \\
  & \quad {}
  + \iO 2\astar \, \p3 \, v
  + \somma j15 \iO \f3j \, \p j \, v
  = 0\,,
  \label{terzaad}
  \\[2mm]
  & - \iO \dt\p4 \, v 
  + \iO \nabla\p4 \cdot \nabla v
  - \chiphi \iO \nabla\p1 \cdot \nabla v
  - \cn \iO \p4 \, v
  = 0\,,
  \label{quartaad}
  \\[2mm]
  & - \< \dt\p5 , v >
  + \iO \nabla\p5 \cdot \nabla v
  - \chia \iO \astar \, \nabla\p3 \cdot \nabla v
  \non
  \\
  & \quad {}
  + \iO \astar \, \juerg{\p5} \, v
  + \somma j15 \iO \f5j \p j \, \, v
  = 0\,,
  \label{quintaad}
\end{align}
\juerg{for} every $v\in V$ and \aet,
and the final conditions
\Beq
  \p1(T) = \b2 \, \bigl( \phistar(T)-\phiO \bigr)
  \aand
  (\p3,\p4,\p5)(T) = (0,0,0) \,.
  \label{cauchyad}
\Eeq
\Accorpa\Pblad primaad cauchyad

\Bthm
Let $\ustar\in\Uad$.
Then, the adjoint problem \Pblad\ has a unique solution satisfying \Regsoluzad.
\Ethm

\Bdim
As for the existence of a solution, also for this problem one can start from a Faedo--Galerkin scheme
constructed by means of the eigenfunctions of the Laplace operator with homogeneous Neumann boundary condition\an{s}.
However, for brevity, we just perform \juerg{the relevant} formal estimates.

\step
First a priori estimate

We test the above equations by $\p1$, $\p1-\p2$, $\p3$, $\p4$, and~$\p5$, respectively.
In addition, we test \eqref{terzaad} by~$-\Delta\p3$.
We obtain~\juerg{the identities}
\begin{align}
  & - \frac 12 \, \ddt \iO |\p1|^2
  + \iO \nabla\p2 \cdot \nabla\p1
  + \somma j15 \iO \f1j \, \p j \, \p1
  = \iO \f10 \, \p1\,,
  \non
  \\
  \separa
  & \iO |\nabla\p1|^2
  - \iO \nabla\p1 \cdot \nabla\p2
  + \iO |\p2|^2
  = \iO \an{\p2 \, \p1\,,}
  \non
  \\
  \separa
  & - \frac 12 \, \ddt \iO |\p3|^2
  + \iO |\nabla\p3|^2
  - \chia \iO (\nabla\sigstar \cdot \nabla\p3) \p3
  + \iO 2\astar \, |\p3|^2
  + \somma j15 \iO \f3j \, \p j \, \p3
  = 0 \,,
  \non
  \\
  \separa
  & - \frac 12 \, \ddt \iO |\p4|^2
  + \iO |\nabla\p4|^2
  - \chiphi \iO \nabla\p1 \cdot \nabla \p4
  - \cn \iO |\p4|^2
  = 0\,,
  \non
  \\
  \separa
  & - \frac 12 \, \ddt \iO |\p5|^2
  + \iO |\nabla\p5|^2
  - \chia \iO \astar \, \nabla\p3 \cdot \nabla\p5
  + \iO \astar \,\juerg{|\p5|^2}
  + \somma j15 \iO \f5j \, \p j \, \p5
  = 0\,,
  \non 
  \\
  \separa
  & - \frac 12 \, \ddt \iO |\nabla\p3|^2
  + \iO |\Delta\p3|^2
  + \chia \iO (\nabla\sigstar \cdot \nabla\p3) \Delta\p3
  \non
  \\
  & \quad {}
  - \iO 2\astar \, \p3 \, \Delta\p3
  - \somma j15 \iO \f3j \, \p j \, \Delta\p3
  = 0 \,.
  \non
\end{align}
Now, we sum up, notice some cancellations, and rearrange a little.
\juerg{Then we obtain the \lhs}
\begin{align}
  & - \frac 12 \, \ddt \iO \bigl(
    |\p1|^2 
    + |\p3|^2
    \an{{} + |\nabla\p3|^2}
    + |\p4|^2
    + |\p5|^2
  \bigr)
  \non
  \\
  & \quad {}
  + \iO \bigl(
  |\nabla\p1|^2
  \an{{} + |\p2|^2}
  + |\nabla\p3|^2
 \an{{} + |\Delta\p3|^2}
 + |\nabla\p4|^2
 	\juerg{+ \astar\,|\p5|^2}
 + |\nabla\p5|^2
  \bigr)\,,
  \label{lhsad}
\end{align}
\juerg{where we recall that $\,\astar\,$ is nonnegative.
We now have} to estimate the terms of the corresponding \rhs,
with the intention of applying the (backward-in-time) Gronwall lemma.
However, since for many of them it suffices to apply Young's inequality, 
we just deal with the ones that need some other treatment.
We consider only one of the easy integrals, using the assumption that $\chiphi\in(0,1)$: \juerg{we have}
\Beq
  \chiphi \iO \nabla\p1 \cdot \nabla \p4
  \leq \frac 12 \iO |\nabla\p1|^2
  + \frac 12 \iO |\nabla\p4|^2\,,
  \non
\Eeq
and the last two integrals are dominated by the corresponding \juerg{expressions occurring in} \eqref{lhsad}.
As for the nontrivial terms, we recall \pier{the regularity for $\nabla\sigstar$ ensured by \eqref{regsig}} and that $\astar\in\LQ4$,
which implies that the function $t\mapsto\norma{\astar(t)}_4$ belongs to $L^4(0,T)$.
We also account for the second \juerg{inequality in} \eqref{lady}.
Hence, we have~that
\begin{align}
  & \chia \iO (\nabla\sigstar \cdot \nabla\p3) \p3
  \leq \norma{\nabla\sigstar}_\infty \, \norma{\nabla\p3}_2 \, \norma{\p3}_2
  \leq c \pier{{}\norma{\sigstar}_{\Hx3} \Bigl(\iO |\nabla\p3|^2
  + \iO |\p3|^2 \Bigr){}},
  \non
  \\[2mm]
  & \chia \iO \astar \, \nabla\p3 \cdot \nabla\p5
  \leq \norma\astar_4 \, \norma{\nabla\p3}_4 \, \norma{\nabla\p5}_2
  \non
  \\
  & \leq \frac 12 \iO |\nabla\p5|^2 
  + c \, \norma\astar_4^2 \, \normaV{\p3} \bigl( \norma{\p3} + \norma{\Delta\p3} \bigr)
  \non
  \\
  & \leq \frac 12 \iO |\nabla\p5|^2 
  + c \, \norma\astar_4^2 \, \normaV{\p3}^2
  + c \, \norma\astar_4^2 \, \normaV{\p3} \, \norma{\Delta\p3} 
  \non
  \\
  & \leq \frac 12 \iO |\nabla\p5|^2 
  + \frac 14 \iO |\Delta\p3|^2
  + c \, \bigl( \norma\astar_4^2 + \norma\astar_4^4 \bigr) \iO (|\p3|^2 + |\nabla\p3|^2) ,
  \non
  \\[2mm] 
  & - \chia \iO (\nabla\sigstar \cdot \nabla\p3) \Delta\p3
  \leq \norma{\nabla\sigstar}_\infty \, \norma{\nabla\p3}_2 \, \norma{\Delta\p3}_2
  \non
  \\
  &  \leq \frac 14 \iO |\Delta\p3|^2
  + c \pier{{}\norma{\sigstar}_{\Hx3}^2{}}\iO |\nabla\p3|^2 \,.
  \non
\end{align}
By combining \an{the above estimates}, \pier{observing that the function $t\mapsto \norma{\sigstar(t)}_{\Hx3}^2$ is bounded in $L^1(0,T)$,} integrating over $(t,T)$ with respect to time, and applying Gronwall's lemma, we conclude~that
\begin{align}
  & \norma{\p1}_{\L\infty H\cap\L2V}
  + \norma{\p2}_{\L2H}
  + \norma{\p3}_{\L\infty V\cap\L2W}
  \non
  \\
  & \quad {}
  + \norma{\p4}_{\L\infty H\cap\L2V}
  + \norma{\p5}_{\L\infty H\cap\L2V}
  \leq c \,,
  \label{1stimaad}
\end{align}
and a comparison in \eqref{secondaad} and in \eqref{terzaad}\an{, along with elliptic regularity,} yields that also
\Beq
  \norma{\p1}_{\L2W}
  + \norma{\p3}_{\H1H}
  \leq c \,.
  \label{da1stimaad}
\Eeq

\step
Second a priori estimate

We test \eqref{primaad} by $\p2$ and \eqref{secondaad} by $-\dt\p1$ and sum up.
Since a cancellation occurs, we obtain~that
\Beq
  \iO |\nabla\p2|^2
  - \frac 12 \, \ddt \iO |\nabla\p1|^2
  = - \somma j15 \iO \f1j \, \p j \, \p2
  + \iO \f10 \, \p2 \,.
  \non
\Eeq
By virtue of the last assumption in \eqref{HPcost}, we notice that $\p1(T)\in V$.
Thus, integrating with respect to time, and applying the Gronwall lemma, we immediately infer that
\Beq
  \norma{\nabla\p1}_{\L\infty H}
  + \norma{\nabla\p2}_{\L2H}
  \leq c \,.
  \non
\Eeq
\juerg{Therefore, we can conclude from \eqref{1stimaad} that}
\Beq
  \norma{\p1}_{\L\infty V}
  + \norma{\p2}_{\L2V}
  \leq c \,.
  \label{2stimaad}
\Eeq

\step
Consequences

Next, we recall that $\p1(T)\in V$ and $\p4(T)=0$.
Hence, we can first compare in \eqref{secondaad}, and then apply the parabolic regularity in \eqref{terzaad} and~\eqref{quartaad},
to see~that
\Beq
  \norma{\p1}_{\L2{\Hx3}}
  + \norma{\p3}_{\H1H}
  + \norma{\p4}_{\H1V\cap\L2{\Hx3}}
  \leq c \,.
  \label{da2stimaad}
\Eeq

\step
Further a priori estimates

Finally, we test \eqref{primaad} and \eqref{quintaad} by a generic $v\in\L2V$ and integrate over~$(0,T)$.
Owing to the previous estimates, we easily conclude~that
\Beq
  \norma{\dt\p1}_{\L2\Vp}
  + \norma{\dt\p5}_{\L2\Vp}
  \leq c \,.
  \label{3stimaad}
\Eeq
We just comment on a term involved in the second test:
we have~that
\Beq
  \chia \intQ \astar \nabla\p3 \cdot \nabla v
  \leq \norma\astar_{\LQ4} \, \norma{\nabla\p3}_{\LQ4} \, \norma{\nabla v}_{\LQ2}
  \leq c \, \norma v_{\L2V} \,.
  \non
\Eeq
This concludes the formal proof of the existence of a solution $\soluzad$ satisfying \Regsoluzad.

\step
Uniqueness

By linearity, we just have to consider the homogeneous problem, i.e.,
we replace $\f10$ by~zero and assume $\p1(T)=0$ in place of the 
first \pier{terminal} value condition in~\eqref{cauchyad}.
Once more, we make one of the formal estimates rigorous.
Namely, we come back to the derivation of \eqref{1stimaad},
where we have tested \eqref{terzaad} by~$-\Delta\p3$.
Instead of doing this, we account for the regularity of the solution,
write \eqref{terzaad} in its strong form,
multiply it by $-\Delta\p3$, and integrate over~$\Omega$.
Then, the same estimates can be performed.
After applying the Gronwall lemma, we conclude that $\soluzad=(0,0,0,0,0)$.
\Edim

Our final result is the first-order necessary optimality condition, expressed in the following theorem\an{.}

\Bthm
\label{Goodnc}
Let $\ustar\in\Uad$ \juerg{be an optimal control}, and let $\soluzad$ be the solution to the associated adjoint problem \Pblad.
Then, there holds the following variational inequality:
\Beq
  \intQ (\p3 + \b3\ustar) (u-\ustar) \geq 0
  \quad \hbox{for every $u\in\Uad$} \,.  
  \label{goodnc}
\Eeq
In particular, \pier{being $\b3>0$,} the optimal control $\ustar$ is the $L^2$-orthogonal projection of $-\p3/\b3$ on~$\Uad$.
\Ethm

\Bdim
We fix $u\in\Uad$, set $h:=u-\ustar$,
and consider both the linearized problem \Pbll\ associated with $\ustar$ and $h$ and the adjoint problem.
We test the equations of the former by $\p1$, $\p2$, $\p3$, $\p4$ and~$\p5$, respectively,
and those of the latter by $-\psi$, $-\eta$, $-\al$, $-\nu$ and~$-\omega$, respectively.
By also recalling \accorpa{deff10}{deff5j}, we have~that
\begin{align}
  & \< \dt\psi , \p1 >
  + \iO \nabla\eta \cdot \nabla\p1
  - \chiphi \iO \nabla\nu \cdot \nabla\p1
  = - m \iO \psi \, \p1
  + \iO \h'(\phistar) \, \psi \, \p1 \,,
  \non
  \\
  \separa
  & \iO \nabla\psi \cdot \nabla\p2
  + \iO F''(\phistar) \, \psi \, \p2 
  = \iO \eta \, \p2 \,,
  \non
  \\
  \separa
  & \< \dt\al , \p3 >
  + \iO \nabla\al \cdot \nabla\p3
  - \chia \iO (\al \nabla\sigstar + \astar \nabla\omega) \cdot \nabla\p3
  = \iO (\al - 2\astar\al + h) \p3\,,
  \non
  \\
  & \iO \dt\nu \, \p4
  + \iO \nabla\nu \cdot \nabla\p4
  - \chiphi \iO \psi \, \p4
  = \iO (\cphi\psi + \cn\nu + \csig\omega) \p4 \,,
  \non
  \\
  \separa
  & \iO \dt\omega \, \p5
  + \iO \nabla\omega \cdot \nabla\p5
  = \iO (-\omega + \chia\al - \al\sigstar - \astar\omega) \p5\,,
  \non
  \end{align}
\pier{and}
   \begin{align}
  & \< \dt\p1 , \psi >
  - \iO \nabla\p2 \cdot \nabla\psi
  - \iO \bigl(
     (m-\h'(\phistar)) \p1
     + F''(\phistar) \, \p2
     - (\chiphi+\cphi) \p4
  \bigr) \psi
  \non
  \\
  & = - \b1 \iO (\phistar-\phiQ) \, \psi\,,
  \non
  \\
  \separa
  & - \iO \nabla\p1 \cdot \nabla\eta
  = - \iO \p2 \, \eta\,,
  \non
  \\
  & \iO \dt\p3 \, \al 
  - \iO \nabla\p3 \cdot \nabla\al
  + \chia \iO (\nabla\sigstar \cdot \nabla\p3) \al
  \non
  \\
  & \quad {}
  - \iO 2\astar \, \p3 \, \al
  - \iO \bigl(
    - \p3
    + (\sigstar-\chia) \p5
  \bigr) \al
  = 0\,,
  \non
  \\
  \separa
  & \iO \dt\p4 \, \nu 
  - \iO \nabla\p4 \cdot \nabla\nu
  + \chiphi \iO \nabla\p1 \cdot \nabla\nu
  + \cn \iO \p4 \, \nu
  = 0\,,
  \non
  \\
  \separa
  & \< \dt\p5 , \omega >
  - \iO \nabla\p5 \cdot \nabla\omega
  + \chia \iO \astar \, \nabla\p3 \cdot \nabla\omega
  - \iO \astar \,\juerg{\p5} \, \omega
  - \iO \bigl(
    - \csig \, \p4 + \p5
  \bigr) \, \omega
  = 0 \,.
  \non
\end{align}
At this point, we take the sum of all these identities.
\juerg{Just a few terms do not cancel out}. 
Namely, we find~that
\begin{align}
  & \< \dt\psi , \p1 >
  + \< \dt\p1 , \psi >
  + \< \dt\al , \p3 >
  + \< \dt\p3 , \al >
  + \< \dt\nu , \p4 >
  + \< \dt\p4 , \nu >
  \non
  \\
  & \quad {}
  + \< \dt\omega , \p5 >
  + \< \dt\p 5, \omega >
  = \iO h \, \p3
  - \b1 \iO (\phistar-\phiQ) \, \psi \,.
  \non
\end{align}
Now, we integrate over $(0,T)$ and apply the well-known integration-by-parts formula
for functions belonging to $\H1\Vp\cap\L2V$.
On account of the initial and final conditions \eqref{cauchyl} and \eqref{cauchyad}, and recalling the choice of~$h$, we then conclude~that
\Beq
  \b2 \iO \bigl( \phistar(T)-\phiO \bigr) \, \psi(T)
  = \intQ (u-\ustar) \, \p3
  - \b1 \juerg{\intQ} (\phistar-\phiQ) \, \psi \,.
  \non
\Eeq
Combining this identity with \eqref{badnc},
we obtain \eqref{goodnc}, and the proof is complete.
\Edim

\juerg{
\Brem
Since the optimal control problem under study is nonconvex, it may have many local minima. We claim
that the variational inequality \eqref{goodnc} has to be valid also for all such locally optimal controls.
In this connection, recall that a control $\ustar\in\Uad$ is termed \emph{locally optimal in the sense
of} $L^p(Q)$ for $\pier{1\leq p\leq\infty}$ if and only if there exists some $\varepsilon>0$ such that
$${\cal J}(\ustar,\calS_1(\ustar))\,\le\,{\cal J}(u,\calS_1(u)) \quad\mbox{for all \,$u\in\Uad$ with $\,\|u-\ustar\|_{L^p(Q)}\le\varepsilon$}.
$$
Note also that every locally optimal control in the sense of $L^p(Q)$ for some $\pier{1\leq p < \infty}$ is also locally optimal
in the sense of $L^\infty(Q)$. Now, it is easily seen that any locally optimal control in the sense of $L^\infty(Q)$ 
satisfies the variational inequality \eqref{badnc}. Hence, by the same argument as in the preceding proof, it must satisfy
also \eqref{goodnc}. 
\Erem
}

\vskip 6mm
\noindent{\bf Acknowledgements}

\noindent
\pier{%
This research was supported by the MIUR-PRIN Grant 2020F3NCPX “Mathematics for industry
4.0 (Math4I4)”. In addition, PC and AS acknowledge their affiliation to the GNAMPA (Gruppo Nazionale per l’Analisi
Matematica, la Probabilità e le loro Applicazioni) of INdAM (Istituto Nazionale di Alta
Matematica). AS has been also supported by the “MUR GRANT Dipartimento di Eccellenza”~2023-2027.}


\End{document}
